%%Version 28.4.11 efter korrekturlaesning

\input amstex
\input amsppt.sty
\magnification1000
\vsize=22truecm
\hsize=15.5truecm
\NoBlackBoxes

\def\Ami{A_{\min}}
\def\Ama{A_{\max}}
\def\inj{\text{\rm i}}
\def\wA{\widetilde A}
\def\Op{\operatorname{Op}}
\def\supp{\operatorname{supp}}
\def\comega{\overline\Omega }
\def\simto{\overset\sim\to\rightarrow}
\def\ang#1{\langle {#1} \rangle}
\def\rnp{{\Bbb R}^n_+}

\def\crpm{\overline{\Bbb R}_\pm}
\def\rp{ \Bbb R_+}
\def\crp{\overline{\Bbb R}_+}
\def\crm{\overline{\Bbb R}_-}

\define\pr{\operatorname{pr}}

\def\d{d\!@!@!@!@!@!{}^{@!@!\text{\rm--}}\!}

\document
\topmatter
\title
The mixed boundary value problem, Krein resolvent formulas and
spectral asymptotic estimates
\endtitle
\author Gerd Grubb \endauthor
\affil
{Department of Mathematical Sciences, Copenhagen University,
Universitetsparken 5, DK-2100 Copenhagen, Denmark.
E-mail {\tt grubb\@math.ku.dk}}\endaffil
\abstract
For a second-order symmetric strongly elliptic
operator $A$ on a smooth bounded open set in ${\Bbb R}^n$, the mixed problem is defined by a Neumann-type condition on a
part $\Sigma _+$ of the boundary and a Dirichlet condition on the
other part $\Sigma _-$. We show a Kre\u\i{}n resolvent formula, where
the difference between its
resolvent and  the Dirichlet resolvent is expressed in terms of
operators acting on Sobolev spaces over $\Sigma _+$. This is used to
obtain a new Weyl-type spectral asymptotics formula for the resolvent
difference (where upper estimates were known before), namely
$s_j j^{2/(n-1)}\to C_{0,+}^{2/(n-1)}$, where $C_{0,+}$ is
proportional to the area of $\Sigma _+$, in the case where $A$ is
principally equal to the Laplacian.

\endabstract
\subjclass 35J25, 35P20, 47G30, 58J40 \endsubjclass

\keywords 
Mixed boundary condition; Zaremba problem; 
resolvent difference; Dirichlet-to-Neumann operator;  Krein resolvent formula; spectral
asymptotics; weak Schatten class; nonstandard
pseudodifferential operator
 \endkeywords

\endtopmatter
\rightheadtext {Mixed problems}

The  mixed boundary value problem for a second-order
strongly elliptic symmetric operator $A$ on a smooth bounded  open set $\Omega
\subset {\Bbb R}^n$ with
boundary $\Sigma $, in case of the Laplacian also called the Zaremba problem,
is defined by a Neumann-type condition
on a part of the boundary $\Sigma _+$ and a Dirichlet condition on the other
part $\Sigma _-$. It does not have the regularity of standard elliptic
boundary problems (the $L_2$-domain is at best in 
$H^{\frac32-\varepsilon }(\Omega )$).
It has been analyzed with regards to 
regularity and mapping properties e.g.\ in Peetre \cite{P61, P63}, Shamir
\cite{S68}, Eskin \cite{E81}, Pryde \cite{P81},
Rempel and Schulze
\cite{RS83}, Simanca \cite{S87}, Harutyunyan and Schulze \cite{HS08}.

We shall here study it from the point of view of extension theory for
elliptic operators.
There has been a recent revival in the interest for connections between
abstract extension theories for operators in Hilbert space (as
initiated by Krein \cite{K47}, Vishik \cite{V52}, Birman \cite{B62},
Grubb \cite{G68} and others) and
interpretations to boundary value problems for partial differential
operators. Cf.\ e.g.\ 
Amrein and Pearson
\cite{AP04}, 
Pankrashkin
\cite{P06}, Behrndt and Langer
\cite{BL07},  Ryzhov \cite{R07},
Brown, Marletta, Naboko and Wood \cite{BMNW08}, 
Alpay and Behrndt \cite{AB09}, Malamud \cite{M10}, based on boundary
triples theory (as developed from the book of Gorbachuk and Gorbachuk
\cite{GG91} and its sources). Other methods are used in the works of
Brown, Grubb and Wood \cite{BGW09}, \cite{G08}, Posilicano and Raimondi
\cite{PR09}, 
Gesztesy and Mitrea
\cite{GM08, GM09, GM11}
(and their references); see also Grubb \cite{G11, G11a, G11b}
and Abels, Grubb and Wood \cite{AGW11}.
One of the interesting aims has been to derive
Kre\u\i{}n resolvent formulas that link the resolvent of a general
operator with the resolvent of a fixed reference operator by
expressing the difference in terms of  operators connected to the boundary. 

For the mixed problem, a Kre\u\i{}n resolvent formula connecting the 
operator to the
Dirichlet  realization was worked out in \cite{P06}, based on boundary
triples theory. A
different formula results from \cite{G68, G74}, see also \cite{BGW09},
Sect.\ 3.2.5. Observations on the connection with the Neumann
realization were given in \cite{M10}. An upper bound for the spectral
behavior of the resolvent difference was shown by Birman in \cite{B62}.
\medskip

In the present paper we shall work out in detail several Kre\u\i{}n resolvent
formulas for the mixed problem. The primary result is a formula where
the difference between the resolvents for the mixed problem and the
Dirichlet problem is expressed explicitly in terms of operators acting over the
subset $\Sigma _+$; this is based on the universal description from
\cite{G68} in terms of operators between closed subspaces of the
nullspace of the maximal operator. In addition, we show some other
explicit formulas related to those of
 \cite{P06}. Mixed problems for $-\Delta $ on creased domains are briefly
 considered, and we establish a Kre\u\i{}n formula for quasi-convex Lipschitz
 domains as defined in \cite{GM11}.

As an application of our primary formula in the smooth case, 
we show how it leads to a
new result giving a Weyl-type spectral
asymptotic estimate for the resolvent difference, with the constant
defined by an integral over $\Sigma _+$; this sharpens considerably
the upper estimates known earlier.  The proof draws on various results for
nonstandard pseudodifferential operators on $\Sigma $.

\head 1. Introduction \endhead

On a bounded smooth open subset $\Omega $ of ${\Bbb R}^n$ with boundary
$\partial\Omega =\Sigma $, consider a second-order symmetric
differential operator with real coefficients in $C^\infty (\comega)$
and an associated sesquilinear form$$
\align
Au&= -{\sum}_{j,k=1}^n \partial _{ j} (a_{jk}(x)
  \partial _{ k} u) + a_0(x)u,\tag1.1\\
a(u,v)&={\sum}_{j,k=1}^n ( a_{jk}
  \partial _{ k} u,\partial _{ j}v) + (a_0u,v).
\tag1.2\endalign
$$ 
$A$ is assumed strongly elliptic, i.e., $\sum_{j,k=1}^n a_{jk}(x)\xi _j\xi _k\ge c_0|\xi
|^2$ for $x\in\Omega $, $\xi \in{\Bbb R}^n$, with $c_0>0$.

Denote $u|_{\Sigma }=\gamma _0u$, and
${\sum}_jn_j\gamma_0(\partial_{ j}u)=\gamma_1u$, where $\vec
n=(n_1,\dots,n_n)$ is the interior unit normal to the boundary. 
Introduce the conormal derivative $\nu $ and a variant $\chi $ (Neumann-type boundary operators) 
$$
\nu u = {\sum}_{j,k=1}^n n_j \gamma_0 (a_{jk}\partial_{ k} u),\quad \chi u=\nu u-b\gamma _0u.
\tag1.3
$$
$\nu $ enters in the ``halfways Green's formula'' (for sufficiently smooth functions)
$$
(Au,v)_{L_2(\Omega )}-a(u,v)=(\nu u,\gamma _0v)_{L_2(\Sigma )} .
\tag1.4
$$

Consider the realizations $A_\gamma $, $A_\nu $, $A_{\chi }$ resp.\
$A_{\chi ,\Sigma _+}$ of $A$ defined via sesquilinear forms to represent the respective boundary conditions$$
\aligned
\gamma _0u&=0 \text{ on }\Sigma ,\text{ the Dirichlet condition},
\\
\nu u&=0 \text{ on }\Sigma ,\text{ the Neumann condition},\\
\chi u&=0 \text{ on }\Sigma ,\text{ a Robin (Neumann-type) condition},\\
\chi u&=0 \text{ on }\Sigma _+,\;\gamma _0u=0\text{ on
}\Sigma \setminus\Sigma _+,\text{ a {\it mixed condition}};
\endaligned\tag1.5
$$
here $b$ is a bounded measurable real function and $\Sigma _+$ is a closed subset of
 $\Sigma $. These realizations are selfadjoint, and by addition of a
 large constant to $a_0$ we can obtain that
 they have positive lower bounds. Their resolvents are compact
 operators. Note that $A_{\chi }$ equals $A_{\nu }$ for $b=0$.

For a compact operator $B$ in a Hilbert space $H$, $s_j(B)$ denotes the $j$-th eigenvalue of $(B^*B)^{\frac12}$ (the
$j$-th s-number or singular value of $B$), counted
with multiplicities. 

Birman showed in \cite{B62}: 
$$
s_j(A_{\chi } ^{-1}-A_\gamma ^{-1})\text{ and }s_j(A_{\chi ,\Sigma _+}
^{-1}-A_\gamma ^{-1}) \text{ are }O(j^{-2/(n-1)})\text{ for }j\to\infty ;
\tag 1.6 
$$
also valid for exterior domains.  
The estimate for $A_{\chi }
^{-1}-A_\gamma ^{-1}$ was later improved to an asymptotic estimate (in
\cite{G74} and \cite{BS80}, the latter including exterior domains):
$$
\lim _{j\to\infty }s _j(A_{\chi }^{-1}-A_\gamma ^{-1})j^{2/(n-1)}=  C_0^{2/(n-1)},\tag1.7
$$
for smooth $b$, where$$
C_0= 
\tfrac1{(n-1)(2\pi )^{n-1}}\int_{\Sigma }\int_{|\xi
'|=1} (\|\tilde k^0\|_{L_2(\rp)}|p^0|^{1/2})^{n-1}
\,d\omega (\xi ') dx';\tag1.8
$$
this has been extended to nonsmooth $b$ in \cite{G11a} (the
ingredients in the formula are explained around  Th.\ 2.4 there).
For  the difference with $A_{\chi ,\Sigma _+}^{-1}$ an asymptotic estimate
does not seem to have been obtained before; it is one of the aims of the
present paper. 

In Section 2, we briefly recall some elements of the old extension theory from
\cite{G68, G74}. In Section 3, we show how the method of Birman
\cite{B62} can be used in combination with later estimates to make a small
improvement of his result for mixed problems, valid for nonsmooth $b$
and $\Sigma _+$. 

In Section 4, 
we analyze the
structure of $A_{\chi ,\Sigma _+}$ in terms of the characterization from
\cite{G68} in more detail,   
describing the operator
$L^\lambda \colon X\to X^*$ that $A_{\chi ,\Sigma _+}-\lambda $
corresponds to when $\lambda \in \varrho (A_\gamma )$ (the resolvent set):

 \proclaim{Theorem A} When $b$ and the subset $\Sigma _+$ are smooth, then
 $X=H^{-\frac12}_0(\Sigma _+)$, and $L^\lambda $ acts like minus the 
Dirichlet-to-Neumann
 pseudodifferential operator truncated to $\Sigma _+$, $-P^\lambda _{\gamma ,\chi
 ,+}=-r^+\chi K^\lambda _\gamma e^+$, with domain $D(L^\lambda )\subset
 H^{1-\varepsilon }_0(\Sigma _+)$ (any $\varepsilon >0$); here $K^\lambda
 _\gamma $ is the Poisson operator for the Dirichlet problem for $A-\lambda $. 

For $\lambda \in \varrho
 (A_{\chi ,\Sigma _+})\cap \varrho (A_\gamma )$ there is a  Kre\u\i{}n resolvent formula:
$$
(A_{\chi ,\Sigma _+}-\lambda ) ^{-1}-(A_\gamma -\lambda ) ^{-1}=-K ^\lambda _{\gamma ,X}(P^\lambda _{\gamma ,\chi
 ,+} )  ^{-1}(K^{\bar\lambda }_{\gamma
,X})^*.\tag1.9%\\
%(A_{\chi ,\Sigma _+} -\lambda )^{-1}-(A_{\nu } -\lambda )  ^{-1}
%=K^\lambda 
%_{{\nu } }
%(P^\lambda  _{{\nu } ,\gamma }-\sigma _+ )^{-1}{K^{\bar\lambda }_{\nu } }^*.\tag1.13
$$
\endproclaim

Several other Kre\u\i{}n-type formulas are shown involving the Poisson operators for
the Dirichlet or Neumann problems.

In Section 5, we restrict the attention to operators principally like
the Laplacian. Here we use methods for
nonstandard pseudodifferential operators to deduce from (1.9):

\proclaim{Theorem B} When $A=-\Delta +a_0(x)$, then for any $\lambda \in \varrho (A_{\chi ,\Sigma _+})\cap
\varrho (A_\gamma )$, 
$$
\lim _{j\to\infty }s _j((A_{\chi ,\Sigma _+}-\lambda )^{-1}-(A_\gamma -\lambda )^{-1})j^{2/(n-1)}=  C_{0,+}^{2/(n-1)},\tag1.10
$$
where $C_{0,+}$ is a constant proportional to the area of $\Sigma _+$; $$
C_{0,+}= 
\tfrac1{(n-1)(2\pi )^{n-1}}\int_{\Sigma _+}\int_{|\xi
'|=1} (\|\tilde k^0\|_{L_2(\rp)}|p^0|^{1/2})^{n-1}
\,d\omega (\xi ') dx'.\tag1.11
$$
\endproclaim

Remark 3.3 and 
Section 4.3 give informations on cases where $\Omega $ is not smooth.

A general technique for extending the estimates to exterior
domains can be found in \cite{G11}.

\head {2. Preliminaries}\endhead 

\subhead 2.1 Definition of the operators \endsubhead

The spaces
$H^s(\Omega )$, $H^s(\Sigma )$ are
the standard Sobolev spaces, with the norm denoted $\|u\|_s$;
$H^s_0(\Omega )$ (or $H^s_0(\comega)$) stands for the space of 
distributions in $H^s({\Bbb
R}^n)$ with support in $\comega$. We use the notation $(\cdot,\cdot)_{-s,s}$ for the sesquilinear duality between
$H^{-s}(\Sigma )$ and $H^{s}(\Sigma )$, $s\in{\Bbb R}$; it reduces to
the $L_2$-scalar product when applied to functions in $L_2(\Sigma )$.

It is known e.g.\ from Lions and Magenes \cite{LM68} that $\gamma _0$ resp.\
$\gamma_1,\nu $ extend to continuous mappings from  
$ H^s(\Omega )\cap D(\Ama)$ to 
$H^{s-\frac12}(\Sigma )$ resp.\ $H^{s-\frac32}(\Sigma )$,
 any $s\ge 0$, allowing extensions of Green's formulas. In particular, for $u\in  H^1(\Omega )\cap
D(\Ama)$, $v\in H^1(\Omega )$,  (1.4) holds with the scalar product
in $L_2(\Sigma )$ replaced by the sesquilinear duality between
$H^{-\frac12}(\Sigma )$ and $H^\frac12(\Sigma )$. 

The realizations of $A$ are the linear operators $\wA$ satisfying
$\Ami\subset \wA\subset \Ama$, where $\Ami$ and $\Ama$ act like $A$
with domains $
D(\Ami)=H^2_0(\Omega )$ resp.\ $D(\Ama)=\{u\in L_2(\Omega )\mid Au\in L_2(\Omega )\}$;
$\Ami$ is the closure of $A|_{C_0^\infty }$, and $\Ama=\Ami^*$.

Our assumptions imply that
$$
a(u,u)\ge c\|u\|^2_{H^1(\Omega )}-k\|u\|^2_{L_2(\Omega )}\text{ for }u\in H^1(\Omega ),\tag2.1
$$
with $c>0$, $k\ge 0$.
Then the realizations $A_\gamma $, etc., can all be defined via variational
constructions from sesquilinear forms, namely:
$$
\aligned
a_\gamma (u,v)&=a(u,v)\text{ on }D(a_\gamma )=H^1_0(\Omega )\text{
leads to }A_\gamma ,\\
a_\nu (u,v)&=a(u,v)\text{ on }D(a_\nu  )=H^1(\Omega )\text{
leads to }A_\nu ,\\
a_{\chi } (u,v)&=a(u,v)+(b\gamma _0u,\gamma _0v)_{L_2(\Sigma )}\text{ on }D(a_{\chi }  )=H^1(\Omega )\text{
leads to }A_{\chi } ,\\
a_{\chi ,\Sigma _+} (u,v)&=a(u,v)+(b\gamma _0u,\gamma _0v)_{L_2(\Sigma
_+)}\text{ on }
 D(a_{\chi ,\Sigma _+}  )=H^1_{\Sigma _+}(\Omega )\text{ leads to }A_{\chi ,\Sigma _+} ;
\endaligned\tag 2.2$$ 
here
$$
H^1_{\Sigma _+}(\Omega ) =\{u\in H^1(\Omega )\mid \supp \gamma _0u\subset
\Sigma _+\}.\tag 2.3
$$
The last case (that covers the two preceding cases when $\Sigma
_+=\Sigma $ or $b=0$) is explained below.
Since 
$\|\gamma _0u\|^2_{L_2(\Sigma )}\le c'\|u\|^2_{\frac34}\le \varepsilon
\|u\|^2_{1}+C(\varepsilon )\|u\|^2_0$ for any $\varepsilon
$, we infer from (2.1) that when $K$ is a constant
$\ge\operatorname{ess} \sup |b(x)|$, $a_{\chi }(u,u)\ge
a(u,u)-K\|\gamma _0u\|^2_0$, and hence 
$$
a_{\chi }(u,u)
\ge c_1\|u\|^2_1-k_1\|u\|^2_0,\text{ for }u\in H^1(\Omega ),
$$
where $c_1<c$ is close to $c$ and $k_1\ge k$ is a large constant.
Then each of the sesquilinear forms in (2.2) satisfies such an
inequality on its domain. Defining $\chi _K=\nu +K\gamma _0$ (the case
$b=-K$), we also have that $a_{\chi _K}(u,v)=a(u,v) -K(\gamma
_0u,\gamma _0v)_\Sigma $ satisfies such an inequality. We can (after a fixed choice of the constant
$K$) replace $A$ by
$A+k_1$, i.e.\ add the constant $k_1$ to the coefficient $a_0$ in (2.1);
then all the resulting sesquilinear forms, including $a_{\chi _K}$,
are positive. {\it For
simplicity, $A+k_1$ and $a(u,v)+k_1\cdot(u,v)$ will in the
following again be
denoted $A$ and $a(u,v)$.}

We now recall the construction of $A_{\chi ,\Sigma _+}$.
The sesquilinear form $a_{\chi ,\Sigma _+}$ on $V=H^1_{\Sigma _+}(\Omega )$ in $H=L_2(\Omega )$
defines an operator $A_{\chi ,\Sigma _+}$ by 
$$
\aligned
D(A_{\chi ,\Sigma _+})&=\{u\in V\mid \exists
f\in H\text{ such that }a_{\chi ,\Sigma _+}(u,v)=(f,v)\text{ for all
}v\in V\},\\
A_{\chi ,\Sigma _+}u&=f.
\endaligned\tag2.4
$$
By J.L.\ Lions' version of the Lax-Milgram
lemma, as recalled e.g.\ in \cite{G09}, Sect.\ 12.4, this defines a
selfadjoint operator with the same lower bound as $a_{\chi ,\Sigma _+}$.
Clearly,
$A_{\chi ,\Sigma _+}$ extends $A|_{C_0^\infty }$, hence $\Ami$, and in view of the
selfadjointness is a restriction of $\Ami^*=\Ama$, so it is a realization of $A$.
By (1.4),
$$
(Au,v)-a_{\chi ,\Sigma _+}(u,v)=(\nu u,\gamma _0v)_{-\frac12,
\frac12}
-(b\gamma _0u,\gamma _0v)_{L_2(\Sigma )}=(\chi u,\gamma _0v)_{-\frac12,\frac12},\tag2.5
$$ when $v\in H^1_{\Sigma _+}(\Omega )$. Thus, when $u\in D(\Ama)\cap H^1_{\Sigma _+}(\Omega )$, $
a_{\chi ,\Sigma _+}(u,v)=(Au,v)$ holds for all $v\in H^1_{\Sigma _+}(\Omega )$
precisely when the distribution $\chi u$ vanishes on
the $H^\frac12$-functions supported in $\Sigma _+$.
In this sense, $A_{\chi ,\Sigma _+}$ represents the boundary condition $\gamma
_0u=0\text{ on }\Sigma \setminus \Sigma
_+$, $\chi u=0\text{ on }\Sigma _+.$ 

The boundary condition can be made more explicit when  
$\Sigma _+$ is a smooth subset of $\Sigma $. We then set
$\Sigma _-=\Sigma \setminus  {\Sigma _+^\circ}$, and have  that
$\Sigma =\Sigma _+\cup \Sigma
_-$, with $\Sigma _+^\circ\cup\Sigma _-^\circ$ dense in $\Sigma $.
Then for $s\in{\Bbb R}$, we denote by
$H^s_0(\Sigma _+)$  the closed subspace of $H^s(\Sigma )$
consisting of the elements with support in $\Sigma _+$. Here $C_0^\infty
(\Sigma _+^\circ)$ is a dense subspace, and it should be noted that
for $s+\frac12\in{\Bbb N}$, the space is different from the space
obtained by closure of  $C_0^\infty (\Sigma _+^\circ)$ in $H^s(\Sigma
_+^\circ)$. For $s\in {\Bbb R}$, the latter space $H^s(\Sigma _+^\circ)$ consists of the restrictions to $\Sigma _+^\circ$ of
distributions in $H^s(\Sigma )$, provided with the quotient norm. 
The spaces $H^s_0(\Sigma _+)$ and  $H^{-s}(\Sigma
_+^\circ)$ are dual with respect to an extension of the $L_2$ scalar product, for all $s\in{\Bbb R}$.  

\proclaim{Lemma 2.1}
When $\Sigma _+$ is smooth,
$$
D(A_{\chi ,\Sigma _+})=
\{u\in H^1(\Omega )\cap D(\Ama)\mid  \gamma _0u\in H^\frac12_0( \Sigma _+),\,
\chi u=0\text{ on }\Sigma _+^\circ\}.\tag2.6
$$
\endproclaim

\demo{Proof} Note first that $\gamma _0H^1_{\Sigma _+}(\Omega
)=H^\frac12_0(\Sigma _+)$, since $\gamma _0H^1(\Omega
)=H^\frac12(\Sigma )$ and $H^\frac12_0(\Sigma _+)$ is the subspace of
$H^\frac12(\Sigma )$ consisting of the functions supported in $\Sigma
_+$. Moreover, $C_0^\infty (\Sigma _+^\circ)$ is dense in
$H^\frac12_0(\Sigma _+)$ and is the image by $\gamma _0$ of the space
of $C^\infty (\comega)$-functions $\psi $ with $\gamma _0\psi $
supported in $\Sigma _+^\circ$.

When $u$ is in the right-hand side of (2.6), then
$$
\ang{\chi u,\gamma _0\psi }=0\text{ for }\gamma
_0\psi \in C_0^\infty (\Sigma _+^\circ);
$$
hence by the denseness of $C_0^\infty (\Sigma _+^\circ)$ in $H^\frac12_0(\Sigma _+)$,
$$
(\chi u,\gamma _0v )_{-\frac12,\frac12}=0\text{ for }v\in H^1_{\Sigma _+}(\Omega ),;
$$
so $u\in D(A_{\chi ,\Sigma _+})$. Conversely, if $u\in D(A_{\chi ,\Sigma _+})$, then $u\in
D(\Ama)\cap H^1_{\Sigma _+}(\Omega )$ implies $\gamma _0u\in
H^\frac12_0(\Sigma _+)$, and since $\chi u$ vanishes on
$H^\frac12$-functions supported in $\Sigma _+$, it vanishes in
particular on $ C_0^\infty (\Sigma _+^\circ)$, i.e., $\nu u-b\gamma
_0u=0$ on $\Sigma _+^\circ$.\qed
\enddemo

\subhead2.2 Abstract extension theories \endsubhead

We shall now connect the operators with the theory of Kre\u\i{}n
\cite{K47}, Vishik \cite{V52}, Birman \cite{B56}, 
Grubb \cite{G68, G70} (the latter also recalled in \cite{BGW09},
the abstract part in \cite{G09}, Ch.\ 13).
The theory of \cite{G68} extends and completes that of \cite{V52} by giving a universal description of  all 
adjoint pairs of extensions
of a dual pair of injective operators. We here just briefly recall how
it describes the extensions $\wA$ of
a symmetric positive operator $\Ami$ with $\Ami\subset \wA\subset
\Ama=\Ami^*$.

The operators act in a Hilbert space $H$ (in the concrete application,
$H=L_2(\Omega )$). Let $A_\gamma $ be the Friedrichs extension of
$\Ami$ (in the application it will be the Dirichlet
realization), and let $Z=\ker \Ama$. Define the decomposition
$$
D(\Ama)=D(A_\gamma )\dot+Z, \text{ with notation }u=u_\gamma +u_\zeta ,\tag2.7$$
where $u_\gamma =\pr_\gamma u=A_\gamma ^{-1}\Ama u$, $u_\zeta
=u-u_\gamma =(1-\pr_\gamma )u=\pr_\zeta u$. This is used in \cite{G68} to show that
there is a 1--1 correspondence between the closed realizations $\wA$
of $A$ and the closed, densely defined operators between closed subspaces of $Z$:
$$
\wA \text{ closed }\longleftrightarrow\cases V,W \subset Z , \text{ closed
subspaces},\\ T \colon V \to W\text{ closed, densely
defined,}\endcases\tag2.8
 $$
 where $D(T)=\pr_\zeta D(\wA)$, $X=\overline{D(T)}$,
 $W=\overline{\pr_\zeta D(\wA^*)}$, and $Tu_\zeta =\pr_W(\Ama u)$
 (here $\pr_W$ denotes orthogonal projection onto $W$). The operator $\wA^* $ corresponds similarly to $T^* \colon W\to V $, and many properties
carry over between $\wA$ and $T$. For example, $\wA $ is invertible
(i.e.\ bijective)
 if and only if $T$ is so, and then we have an
abstract resolvent formula:
$$\wA ^{-1}=A_\gamma 
^{-1}+\inj_{V }T^{-1}\pr_{W},\tag2.9$$
where $\inj_V$ denotes the injection $V\hookrightarrow H$.

In particular, $\wA$ is selfadjoint if and only if:  $V=W$ and 
 $T\colon V\to V$ is selfadjoint. Then in the invertible case,
$$
\wA^{-1}=A_\gamma ^{-1}+\inj_V T^{-1}\pr_V.\tag2.10
$$
Positivity of $\wA$ holds if and only if $T$ is
positive. 

For the positive selfadjoint operators, there is also a connection between the associated sesquilinear
forms. (When $S$ is a positive selfadjoint operator in a Hilbert space
$H$, the associated sesquilinear form $s$ has as its domain $D(s)$ the
completion of $D(S)$ in the norm $(Su,u)^\frac12$, stronger than the $H$-norm; here
$D(s)\subset H$, and the form $s(u,v)$ is the extension by continuity
of $(Su,v)$ to $D(s)$. Then $S$ is defined from $s$ by the Lax-Milgram
construction.) When $\wA$ is positive selfadjoint, corresponding to
the positive selfadjoint operator $T$ in $V$,
the associated sesquilinear form $\widetilde a$ can be written $$
\widetilde a(u,v)=a_\gamma (u_\gamma ,v_\gamma )+t(u_\zeta ,v_\zeta
)\text{ on }D(\widetilde a)=D(a_\gamma )\dot+ D(t),
\tag2.11$$
where $t$ on $D(t)\subset V$ is the sesquilinear form associated with
$T$; the decomposition $u=u_\gamma +u_\zeta $ used here is a
continuous  extension to
$D(a_\gamma )\dot+Z$ of the decomposition (2.7) above.

The description of selfadjoint extensions in terms of sesquilinear
forms is already found in
\cite{K47} and \cite{B56}; \cite{G70} moreover treats nonselfadjoint
extensions. 

Much of the theory holds unchanged if we replace the ``reference
operator''  $A_\gamma $ by
another selfadjoint positive realization of $A$, say $A_\nu $ (which
will in the application be taken as the Neumann realization $A_\nu $). There
is again a decomposition $$
D(\Ama)=D(A_\nu )\dot+Z, \text{ say with notation }u=u_\nu +u_{\zeta,1 },$$
where $u_\nu =\pr_\nu u=A_\nu ^{-1}\Ama u$, $u_{\zeta,1}
=u-u_\nu    =(1-\pr_\nu    )u=\pr_{\zeta ,1}u$, and there is a 1--1 correspondence
$$
\wA \text{ closed }\longleftrightarrow\cases V_1,W_1 \subset Z , \text{ closed subspaces},\\ T_1 \colon V _1\to W_1\text{ closed, densely defined,}\endcases
\tag2.12 $$
 where $D(T_1)=\pr_{\zeta,1} D(\wA)$, $X_1=\overline{D(T_1)}$,
 $W_1=\overline{\pr_{\zeta,1} D(\wA^*)}$, and $T_1u_{\zeta,1}
 =\pr_{W_1}(\Ama u)$;
again $\wA$ is selfadjoint or invertible if and only if $T_1$ is so,
 and in the invertible case,
$$\wA ^{-1}=A_\nu
^{-1}+\inj_{V _1}T_1^{-1}\pr_{W_1}.\tag2.13$$
However, positivity does not in general carry over between $\wA$ and $T_1$, and the information on associated sesquilinear forms does not generalize to
this situation, since those facts depended on $A_\gamma $ being the Friedrichs
extension of $\Ami$.

\subhead 2.3 Concrete boundary conditions. Dirichlet reference operator \endsubhead

We now explain the interpretation to concrete boundary conditions
worked out in \cite{G68, G74}.
Along with (1.4) we have the full Green's formula
$$
(Au,v)_{L_2(\Omega )}-(u,Av)_{L_2(\Omega )}=(\nu u,\gamma _0v)_{L_2(\Sigma
)}-(\gamma _0u,\nu v)_{L_2(\Sigma )}, \text{ for }u,v\in H^2(\Omega );\tag2.14
$$
it extends e.g.\ to $u\in D(\Ama)$, $v\in H^2(\Omega )$, with the $L_2(\Sigma
)$-scalar products replaced by suitable Sobolev space dualities, but it
cannot be extended to $u,v\in D(\Ama)$.

 Denote by $K _\gamma $ resp.\ $K _\nu $ the
Poisson operator solving the Dirichlet problem resp.\ Neumann problem
$$
Au=0\text{ in }\Omega , \text{ with }\gamma _0u=\varphi
,\text{ resp. }\nu u=\psi ;
$$
they have the mapping properties
$$
K _\gamma \colon H^{s-\frac12}(\Sigma )\to H^s(\Omega ),\;
K _\nu \colon H^{s-\frac32}(\Sigma )\to H^s(\Omega ), \text{
 for all }s\in{\Bbb  R}.
$$
In particular, $\gamma _0$ and $\nu $ define {\it homeomorphisms} of
$Z$ onto $H^{-\frac12}(\Sigma )$ resp.\ $H^{-\frac32}(\Sigma )$, with
$K_\gamma $  resp.\ $K_\nu $ acting as inverses. 

Let $\wA$ correspond to $T\colon V\to W$ as in (2.8). Let $X=\gamma
_0(V)$, $Y=\gamma _0(W)$, closed subspaces of $H^{-\frac12}(\Sigma )$,
and introduce the notation for the connecting homeomorphisms
$$
\gamma _{V }\colon
V \simto X,\quad\gamma _{W}\colon
W\simto Y.\tag2.15
$$
By use of these homeomorphisms, $T \colon V \to
W$ is carried over to a map $L \colon X\to Y^*$:
$$
\CD
V      @> \sim >  \gamma _{V }  >   X\\
@VT  VV           @VV  L   V\\
   W  @<  \sim< \gamma _{W}^*<  Y^*\endCD \hskip1cm
   D(L )=\gamma _0 D(T )=\gamma _0 D(\wA).
 $$
In other words, 
$$
L=(\gamma _{W}^*)^{-1}T \gamma _{V }^{-1}.
$$

In the case where $\wA$ is invertible, the abstract resolvent formula (2.9) carries over to
the formula:
$$\wA^{-1}
=A_\gamma ^{-1}+K _{\gamma ,X} 
L^{-1}(K_{\gamma ,Y})^*\tag2.16
$$
where $$
K _{\gamma
  ,X}= \inj_{V  }\gamma _{V }^{-1}\colon X\to
V \subset H,\quad(K _{\gamma
  ,Y})^*= (\gamma _{W }^*)^{-1}\pr_W\colon H\to Y^*;\tag2.17
$$ (2.16) is a {\it Kre\u\i{}n resolvent formula}. In particular, if $V=W=Z$, then $X=Y=H^{-\frac12}(\Sigma
)$, and (2.16) takes the form
$$\wA^{-1}
=A_\gamma ^{-1}+K _{\gamma } 
L^{-1}{K_{\gamma }}^*,\tag2.18
$$
where $L$ goes from $D(L)\subset H^{-\frac12}(\Sigma )$ to
$H^{\frac12}(\Sigma )$.

To see how $L$ enters in a concrete boundary condition for $\wA$ we
define some additional operators, namely
the {\it Dirichlet-to-Neumann} and {\it 
 Neumann-to-Dirichlet} pseudodifferential operators ($\psi $do's)
$P_{\gamma ,\nu }$ and $P_{\nu ,\gamma }$, and
 the associated {\it reduced trace operators} $\Gamma _\nu $ and $\Gamma _\gamma $:
$$
\aligned
P_{\gamma ,\nu }&=\nu K _\gamma ,\text{ $\psi $do of
order 1,} \quad \Gamma
 _\nu = \nu -P _{\gamma ,\nu }\gamma _0\colon
D(\Ama)\to H^\frac12(\Sigma );\\  
P_{\nu ,\gamma  }&=\gamma _0 K_\nu , \text{ $\psi
$do of order $-1$,}\quad \Gamma
 _\gamma = \gamma _0 -P _{\nu,\gamma  }\nu \colon
D(\Ama)\to H^\frac32(\Sigma ).
\endaligned\tag2.19$$
(We here use the notation of the pseudodifferential boundary operator
calculus, initiated by Boutet de Monvel \cite{B71} and further
developed in \cite{G84, G96}, see also \cite{G09}.) More generally, $P_{\beta ,\beta '}$ denotes the mapping from $\beta
u$ to $\beta 'u$, when $u\in Z$ is uniquely determined from $\beta u$.

The reduced trace operators are used to establish generalized Green's formulas valid for $u,v\in D(\Ama)$:
$$\aligned
(Au,v)_{L_2(\Omega )}-(u,Av)_{L_2(\Omega )}&=(\Gamma  _\nu u,\gamma _0v)_{\frac12,-\frac12}-(\gamma _0u,\Gamma _\nu v)_{-\frac12,\frac12},\\
(Au,v)_{L_2(\Omega )}-(u,Av)_{L_2(\Omega
)}&=
(\nu u,\Gamma _\gamma v)_{-\frac32,\frac32}-(\Gamma  _\gamma u,\nu v)_{\frac32,-\frac32}.
\endaligned\tag2.20$$
One can then show:

$D(\wA)$ consists of the functions  $u\in D(\Ama)$
that satisfy:
$$
\gamma _0u\in D(L ),\quad (\Gamma _\nu  u,\varphi )_{\frac12,
-\frac12}=(L \gamma _0u,\varphi  )_{Y^*,Y}\text{ for all }\varphi  \in Y.
\tag2.21$$
The second condition may be rewritten as $\inj_Y^*\Gamma _\nu 
u=L \gamma _0u$, where $\inj_Y^*\colon H^\frac12(\Sigma )\to Y^*$ is
the adjoint of $\inj_Y\colon Y\hookrightarrow H^{-\frac12}(\Sigma
)$. By the definition of $\Gamma _\nu $, this can be written:
$$
\inj_Y^*\nu  u= (L +\inj_Y^*P _{\gamma ,\nu })\gamma _0u.\tag2.22
$$
In the case where $X=Y=H^{-\frac12}(\Sigma )$, this is simply a
Neumann-type condition
$$
\nu  u= C\gamma _0u,\text{ where }C=L +P _{\gamma ,\nu }.
\tag2.23$$
In the present paper we are more interested in a genuine subspace case, where
$X=H^{-\frac12}_0(\Sigma _+)$; we return to that below.

\subhead 2.4 Neumann reference operator \endsubhead

For the abstract theory using $A_\nu $ as the reference operator, we
get slightly different but analogous formulas:

Let $\wA$ correspond to $T_1\colon V_1\to W_1$ as in (2.12). We now set $X_1=\nu
(V_1)$, $Y_1=\nu (W_1)$, closed subspaces of $H^{-\frac32}(\Sigma )$, and denote the connecting homeomorphisms
$$
\nu _{V _1}\colon
V_1 \simto X_1,\quad\nu _{W_1}\colon
W_1\simto Y_1.\tag2.24
$$
Now $T_1 \colon V_1 \to
W_1$ is carried over to the map $L_1 \colon X_1\to Y_1^*$ defined by
$$
L_1=(\nu _{W_1}^*)^{-1}T_1 \nu _{V_1 }^{-1}.
\tag2.25$$

In the invertible case, the abstract resolvent formula (2.13) carries over to
the formula:
$$\wA^{-1}
=A_\nu ^{-1}+K _{\nu ,X_1} 
L_1^{-1}{K_{\nu ,Y_1}}^*\tag2.26
$$
where $K _{\nu
  ,X_1}= \inj_{V _1 }\nu _{V_1 }^{-1}\colon X_1\to
V_1 \subset H$, $(K _{\nu
  ,Y_1})^*= (\nu _{W _1}^*)^{-1}\pr_{W_1}\colon H\to Y_1$; another  Kre\u\i{}n resolvent formula. In particular, if $V_1=W_1=Z$, then $X_1=Y_1=H^{-\frac32}(\Sigma
)$, and (2.26) takes the form
$$\wA^{-1}
=A_\nu ^{-1}+K _{\nu } 
L_1^{-1}{K_{\nu }}^*,\tag2.27
$$
where $L_1$ goes from $D(L_1)\subset H^{-\frac32}(\Sigma )$ to
$H^{\frac32}(\Sigma )$.

The interpretation of $\wA$ as defined by a boundary condition is here
based on the second line of (2.20) and goes as follows: $D(\wA)$ consists of the functions  $u\in D(\Ama)$
that satisfy the boundary condition
$$
\nu u\in D(L _1),\quad -(\Gamma _\gamma  u,\varphi )_{\frac32,
-\frac32}=(L _1\nu u,\varphi  )_{Y_1^*,Y_1}\text{ for all }\varphi  \in Y_1.
\tag2.28$$
Here the second condition is rewritten as $\inj_{Y_1}^*\Gamma _\gamma 
u=-L _1\nu u$, or
$$
\inj_{Y_1}^*\gamma _0  u= (-L _1+\inj_{Y_1}^*P _{\nu ,\gamma })\nu u.
\tag2.29$$
In the case where $X_1=Y_1=H^{-\frac32}(\Sigma )$, this is a
``Dirichlet-type'' condition
$$
\gamma _0  u= C_1\nu u,\text{ where }C_1=-L_1 +P _{\nu ,\gamma  }.
\tag2.30$$
We shall see later that the mixed problem can be written in this form (after
a replacement of $\nu $ by $\nu +K\gamma _0$, if necessary).

In the above analysis we assumed $A_\gamma $ resp.\ $A_\nu $ positive,
so that $0\in \varrho (A_\gamma )$ resp.\ $0\in \varrho (A_\nu
)$. Clearly, by addition of real constants to $A$ this covers the
realizations of $A-\lambda $ for $-\lambda $ large positive. The formulation was
just chosen for simplicity of notation; the theory of \cite{G68} in fact works for any
$\lambda \in \varrho (A_\gamma )$ resp.\ $\lambda \in \varrho (A_\nu )
$. For general $\lambda $ one uses the nullspaces $Z_\lambda =\ker(\Ama-\lambda
)$ and $Z_{\bar\lambda }=\ker(\Ama-\bar\lambda )$. For the various spaces, operators
and auxiliary Poisson,
pseudodifferential and trace operators, the $\lambda $-dependence is
indicated by
$$
V_\lambda ,\, W_{\bar\lambda },\, L^\lambda ,\, K^\lambda _{\gamma },\, K^{
\bar\lambda }_\gamma ,\, P^\lambda _{\gamma ,\nu },\, P^\lambda _{\nu
,\gamma },\, \Gamma ^\lambda _\nu , \text{ etc.}
\tag2.31$$
The $\lambda $-dependent formulas are explained in detail in 
\cite{BGW09} (based on methods from \cite{G74}), see also \cite{AGW11}
for notation. 
There is an important point here, namely that $X=\gamma _0V_\lambda $ and
$Y=\gamma _0W_{\bar\lambda }$ are {\it independent of $\lambda $}.
Moreover $D(L^\lambda )=D(L^0)$, and $L^\lambda -L^0$
acts as the {\it bounded} operator $\inj_Y^*(P^0_{\gamma ,\nu
}-P^\lambda _{\gamma ,\nu })$. Related statements hold for
$L_1^\lambda \colon X_1\to Y_1$. The Kre\u\i{}n resolvent formulas have
the form:
$$
\gathered
(\wA-\lambda )^{-1}
=(A_\gamma -\lambda )^{-1}+K ^\lambda _{\gamma ,X} 
(L^\lambda )^{-1}(K^{\bar\lambda }_{\gamma ,Y})^*\text{ when }\lambda
\in \varrho (A_\gamma )\cap \varrho (\wA),\\
(\wA-\lambda )^{-1}
=(A_\nu -\lambda )^{-1}+K ^\lambda _{\nu ,X_1} 
(L_1^\lambda )^{-1}(K^{\bar\lambda }_{\nu ,Y_1})^*\text{ when }\lambda
\in \varrho (A_\nu )\cap \varrho (\wA).
\endgathered\tag2.32
$$

Other Kre\u\i{}n resolvent formulas have been established e.g.\ in
Malamud and Mogilevskii \cite{MM02}, \cite{M10}, Pankrashkin
\cite{P06}, Behrndt and Langer
\cite{BL07},  Alpay and Behrndt \cite{AB09},
Gesztesy and Mitrea \cite{GM08, GM09, GM11}, Brown, Marletta, Naboko
and Wood \cite{BMNW08},
Posilicano  and Raimondi \cite{PR09}.

\example{Remark 2.2} The theory recalled above has, in the study of
``pure'' boundary conditions (of Neumann-type $\nu 
u=C\gamma _0u$ or of Dirichlet-type $\gamma _0u=C_1\nu u$), much in
common with the representations of boundary value problems based on
boundary triples theory. 
It is when
{\it subspaces} $V,W$ of $Z$ occur that our theory differs markedly from the
others, which obtain a generalization by allowing {\it relations} instead of
operators.
\endexample

\head 3. Birman's method revisited\endhead

The  correspondence (2.8) with $A_\gamma $ as reference
operator is used here.
We have that $D(A_\gamma )=H^1_0(\Omega )\cap H^2(\Omega )$ and
$D(a_\gamma )=H^1_0(\Omega )$. For $A_{\chi }$, the decomposition in (2.11)
gives  $D(a_{\chi })=H^1(\Omega )=H^1_0(\Omega )\dot+ Z^1$, where  
$Z^1=Z\cap H^1(\Omega )$. The
corresponding operator $T_{\chi }$ is defined from the
sesquilinear form $t_{\chi }$ obtained by restricting $a_{\chi }$ to  $Z^1$ in $Z$; $T_{\chi }$ is a selfadjoint unbounded
positive operator in $Z$ with domain dense in $Z^1$. For the mixed problem,
$D(a_{\chi ,\Sigma _+})=H^1_0(\Omega )\dot+ Z^1_{\Sigma _+}$, where $Z^1_{\Sigma _+}=Z\cap H^1_{\Sigma _+}(\Omega
)$ (cf.\ (2.3)); the corresponding operator $T_{\chi ,\Sigma _+}$ is a selfadjoint
operator in $Z_{\Sigma _+}=\overline {Z^1_{\Sigma _+}}$ ($L_2(\Omega )$-closure) with domain
dense in
$Z^1_{\Sigma _+}$.

There are bounded, in fact compact, inverses $T_{\chi }^{-1}$ on $Z$, resp.\ $T_{\chi ,\Sigma _+}^{-1}$ on $Z_{\Sigma _+}$.

When a general $T$ is derived from the form $t=\widetilde a|_{D(t)}$ and $T^{-1}$ is compact nonnegative, then the
eigenvalues are determined by the minimum-maximum principle from Rayleigh
quotients:
$$
\mu _j(T ^{-1})=\min_{U\subset D(t
), \dim U=j-1}\;\max_{z\perp U, z\in D( t )\setminus
\{0\}}\;\frac{\|z\|_0^2}{\widetilde a (z,z)}.\tag3.1
$$
This principle was used in Birman \cite{B62} to reduce the proof of upper
estimates of the $\mu _j(T^{-1})$ for each of the boundary conditions (1.5) to simpler cases where it could be found by computation.

We shall here show how the principle leads to a lim sup estimate for the
mixed problem.
Consider $a_{\chi ,\Sigma _+}$ and the Robin case $a_{\chi _K}$ (where
$b$ is replaced by $-K$, cf.\ Section 2.1). 
Let the corresponding operators and forms defined on subspaces of $Z$ be denoted
$T_{\chi ,\Sigma _+}$ and $T_{\chi _K}$, resp.\ 
$t_{\chi ,\Sigma _+}$ and $t_{\chi _K}$. Here $ D(t_{\chi _K})=
 Z^1$, and $D(t_{\chi ,\Sigma _+})=Z^1_{\Sigma _+}\subset Z^1$.  Then
$$\aligned
&\mu _j(T _{\chi ,\Sigma _+}^{-1})=\min_{U\subset Z^1_{\Sigma _+}
, \dim U=j-1}\;\max_{z\perp U, z\in Z^1_{\Sigma _+} \setminus
\{0\}}\;\frac{\|z\|_0^2}{a (z,z)+(b\gamma _0z,\gamma _0z)_{\Sigma
_+}}\\
&\le\min_{U\subset Z^1, \dim U=j-1}\;\max_{z\perp U, z\in Z^1\setminus
\{0\}}\;\frac{\|z\|_0^2}{a (z,z)+(b\gamma _0z,\gamma _0z)_{\Sigma _+}}
\\
&\le\min_{U\subset Z^1, \dim U=j-1}\;\max_{z\perp U, z\in Z^1\setminus
\{0\}}\;\frac{\|z\|_0^2}{a (z,z)-K\|\gamma _0z\|^2_{L_2(\Sigma )}}=\mu _j(T^{-1}_{\chi _K}).
\endaligned\tag3.2
$$

Birman showed in \cite{B62} that $\mu _j(T_{\chi _K}^{-1})$, and hence
also $\mu _j(T_{\chi ,\Sigma _+}^{-1})$, is $O(j^{-2/(n-1)})$ for $j\to\infty $.
It is noteworthy that this included the mixed problem. 

In the finer asymptotic estimate (1.7)--(1.8), $p^0(x',\xi ')$ denotes the principal symbol
of $P_{\gamma ,\nu }$ and $\tilde k^0(x',\xi ',\xi _n)$ is the principal
symbol-kernel of $K_\gamma $; the derivation of the formula is
explained in \cite{G11a}, Th.\ 2.4. Applying (1.7)--(1.8) to $T_{\chi _K}$ we can now
get a lim sup estimate using (3.2): 

\proclaim{Proposition 3.1} The nonzero eigenvalues of
$A_{\chi ,\Sigma _+}^{-1}-A_\gamma ^{-1}$ satisfy, with $C_0$ from {\rm (1.8)},
$$
{\lim \sup}_{j\to\infty }\mu _j(A_{\chi ,\Sigma _+}^{-1}-A_\gamma ^{-1})j^{2/(n-1)}\le   C_0^{2/(n-1)}.\tag3.3
$$
\endproclaim

\demo{Proof} From (1.7) with $b=-K$ follows in view of 
(3.2): 
$$
\aligned
{\lim \sup}_{j\to\infty }\mu _j(T_{\chi ,\Sigma _+}^{-1})j^{2/(n-1)}&\le
{\lim \sup}_{j\to\infty }\mu _j(T_{\chi _K}^{-1})j^{2/(n-1)}\\
&={\lim }_{j\to\infty }\mu  _j(A_{\chi _K}^{-1}-A_\gamma ^{-1})j^{2/(n-1)}
=
C_0^{2/(n-1)};
\endaligned\tag3.4
$$
we have here applied formula (2.10) with $\wA=A_{\chi _K}$. Similarly, $A_{\chi ,\Sigma _+}^{-1}-A_\gamma
^{-1}$ and $T_{\chi ,\Sigma _+}^{-1}$ have the same nonzero eigenvalues, so the
result follows.\qed
\enddemo

We also get a spectral estimate for the eigenvalues of $A_{\chi ,\Sigma _+}^{-1}$
 itself:

\proclaim{Corollary 3.2} The eigenvalues of $A_{\chi ,\Sigma _+}$ satisfy:
$$
\mu _j(A_{\chi ,\Sigma _+} ^{-1})-C_A^{2/n}j^{-2/n}\text{ is
}O(j^{-(1+1/(n+1))2/n})\text{ for }j\to\infty ,\tag3.5
$$
where
$$
C_A=(2\pi )^{-n}\int_{x\in \Omega ,\,a^0(x,\xi )<1}\,dxd\xi .
\tag3.6
$$
\endproclaim 

\demo{Proof} It is known (cf.\ e.g.\ \cite{H85}, Sect. 29.3) that the spectrum of $A_\gamma $ satisfies the asymptotic
estimate
$$
\mu _j(A_\gamma  ^{-1})-C_A^{2/n}j^{-2/n}\text{ is
}O(j^{-3/n})\text{ for }j\to\infty ,\tag3.7
$$
with $C_A$ defined by (3.6) (the spectral estimate is formulated for
the counting function in \cite{H85}, but carries over to the above
form, cf.\ e.g.\ \cite{G96}, Lemma A.5). We shall apply a perturbation result to
this estimate, using (3.3) and (2.10) with $\wA=A_{\chi ,\Sigma _+}$.

Recall from \cite{G84}, Prop.\ 6.1 (or \cite{G96}, Lemma A.6), that
when
$B$ and $B'$ are compact operators satisfying for $j\to\infty$, with
$p>q>0$, $p>r>0$, $c_0\ge 0$,
 $$
s_j(B)-c_0^{1/p}j^{-1/p}\text{ is }O(j^{-1/q}),\quad
s_j(B')\text{ is }O(j^{-1/r}),\tag 3.8
$$
then $B+B'$ satisfies
 $$s_j(B+B')-c_0^{1/p}j^{-1/p}\text{ is }O(j^{-1/q'}),\text{ with }q'=\max\left\{q,p\tfrac{r+1}{p+1}\right\}.\tag 3.9
 $$
We apply the result here with $B=A_\gamma ^{-1}$ and
$B'=A_{\chi ,\Sigma _+}^{-1}-A_\gamma ^{-1}$, so that $p=n/2$, $q=n/3$,
$r=(n-1)/2$. This gives 
$$
q'=\max\big\{\frac{n}3,
\frac n2\cdot\frac{\frac{n-1}2+1}{\frac n2+1}\}=\frac n2\cdot\frac{n+1}{n+2};
$$
 here $1/q'=2/n\cdot (n+2)/(n+1)=(1+1/(n+1))2/n$.\qed
\enddemo

Note that these results hold when $b\in L_\infty (\Sigma )$ and $\Sigma _+$ is
any closed subset of $\Sigma $.

\example{Remark 3.3} Concerning nonsmooth choices of $\Omega $, let us
mention that the basic hypothesis  of Birman in \cite{B62} is that
$\Sigma $ is piecewise $C^2$. This allows edges or creases, cf.\ Section 4.3
below.  Moreover, in the recent translation to English of that
historical paper, the translator M.\ Solomyak states in a supplementing 
comment to Section 2.1 on page 50 that 
the result is valid for Lipschitz domains (at least when $n\ge 3$; the
reservation for $n=2$ seems to be concerned with cutoffs in exterior domains).
\endexample

\head 4. Kre\u\i{}n resolvent formulas for the mixed problem\endhead

\subhead 4.1 A formula relative to the Dirichlet problem \endsubhead

We assume in Sections 4.1, 4.2 and 5 that $\Sigma _+$ is smooth.
First we show a Kre\u\i{}n resolvent formula for $A_{\chi ,\Sigma _+}$
linked with $A_\gamma $.
For simplicity of notation, we do the main calculations in the case
$\lambda =0$ (where the indexation by $\lambda $ is left out); then at
the end we account for the consequences in situations with other
values of $\lambda $.

Recall from Section 2.3 that in the analysis 
with $A_\gamma $ as the reference operator,
$A_{\chi ,\Sigma _+}  $ corresponds to $L \colon X\to X^*$, where
$D(L )=\gamma _0D(A_{\chi ,\Sigma _+})$ and $X$ is its closure in
$H^{-\frac12}(\Sigma )$.
It is seen from (2.6) that $D(L )$ is a subset of
$H^{\frac12}_0(\Sigma _+)$, and it  contains $C_0^\infty (\Sigma
_+^\circ)$ in view of the surjectiveness of $\{\gamma _0,\nu
\}$ from $H^2(\Omega )$ to $H^\frac32(\Sigma )\times
H^{\frac12}(\Sigma )$. Then in fact its closure $X$ in
$H^{-\frac12}(\Sigma )$ satisfies 
$$
X=H^{-\frac12}_0(\Sigma _+),\text{ and hence }X^*=H^\frac12(\Sigma _+^\circ).\tag4.1
$$
We note that the injection $\inj_X\colon X\hookrightarrow H^{-\frac12}(\Sigma
)$ and its adjoint satisfy:
$$
\inj_X=e_{\Sigma _+^\circ}\colon H^{-\frac12}_0(\Sigma _+)\hookrightarrow
H^{-\frac12}(\Sigma ),\quad (\inj_X)^*=r_{\Sigma _+^\circ}\colon
H^\frac12(\Sigma )\to H^\frac12(\Sigma _+^\circ),
$$
where $e_{\Sigma _+^\circ}$ is a well-defined extension of  the
operator that extends functions on $\Sigma ^\circ_+$ by zero on
$\Sigma _-$, and $r_{\Sigma _+^\circ}$ denotes restriction to 
$\Sigma _+^\circ$. We denote $e_{\Sigma _+^\circ}=e^+$ and 
$r_{\Sigma _+^\circ}=r^+$ for short. Since $A_{\chi ,\Sigma _+}$ is bijective, so
is $L$, from $D(L)$ to $H^\frac12(\Sigma _+^\circ)$.

When $u\in D(A_{\chi ,\Sigma _+})$, we see from (2.6), (1.3) that
$\nu u$ equals $b\gamma _0u$ on $\Sigma _+^\circ$ in the distribution
sense, hence since $\Gamma _\nu =\nu -P_{\gamma ,\nu }\gamma _0$, $$
\langle\Gamma _\nu  u,\bar\zeta \rangle=\langle(b-P _{\gamma ,\nu
})\gamma _0u,\bar\zeta \rangle \text{ for }\zeta \in C_0^\infty
(\Sigma _+^\circ).\tag4.2
$$
Since $\gamma _0u\in H^\frac12(\Sigma  )$, which is mapped to
$H^{-\frac12}(\Sigma )$ by $P_{\gamma ,\nu }$, and multiplication by
$b$ preserves  $L_2(\Sigma )$, we have that  $(b-P _{\gamma ,\nu
})\gamma _0u\in H^{-\frac12}(\Sigma )$.

The operator $L $ satisfies, by (2.21),
$$
(L \gamma _0u,\varphi )_{X^*,X}=(\Gamma _\nu  u,\varphi
)_{\frac12,-\frac12}\text{ for all }\varphi \in X; 
$$
in particular, when (4.1) and (4.2) are taken
into account,
$$
(L \gamma _0u,\zeta )_{H^\frac12(\Sigma _+^\circ),H^{-\frac12}_0(\Sigma _+)}=\langle(b-P _{\gamma ,\nu
})\gamma _0u,\bar\zeta \rangle,\text{ for }\zeta \in C_0^\infty
(\Sigma _+^\circ),
$$
so
$$
L \gamma _0u=r^+ (b-P _{\gamma ,\nu
})\gamma _0u,\text{ for }u\in D(A_{\chi ,\Sigma _+}).
$$
Thus $L $ acts as 
$$
L \varphi = r^+ (b-P _{\gamma ,\nu
})e^+\varphi , \text{ for }\varphi \in D(L).\tag4.3
$$

This shows the form of $L $. We need deeper theories to say
more about the 
domain. Here we shall use the study of mixed problems in
Shamir \cite{S68}; in Section 5 we also use Eskin \cite{E81}.
Some smoothness is needed for this; for convenience we take
$b\in C^\infty (\Sigma )$.

\proclaim{Proposition 4.1} When $\Sigma _+$ is smooth, the operator $L$ acts
as in {\rm (4.3)}. When also $b$ is smooth,
it satisfies$$
D(L)\subset H^{1-\varepsilon }_0(\Sigma _+),\text{ any }\varepsilon >0,\tag4.4
$$
and $L^{-1}$ maps $H^\frac12(\Sigma _+^\circ)$ into $H^{1-\varepsilon
}_0(\Sigma _+)$.
\endproclaim

\demo{Proof}
We see  from \cite{S68} that  $D(A_{\chi ,\Sigma
_+})\subset H^{\frac32-\varepsilon
}(\Omega )$, as follows: First Shamir shows this in  Th.\ 3.1  of \cite{S68}
for
the constant-coefficient case of $-\Delta +\alpha ^2$ on a half-space
with mixed Dirichlet and Neumann boundary conditions. Subsequently the
statement is
extended to variable coefficients and bounded domains in the proof 
of Lemma 5.1 of \cite{S68} (when we recall that the domain is a priori
contained in $H^1(\Omega )$). Since $\gamma _0H^{\frac32-\varepsilon
}(\Omega )=H^{1-\varepsilon }(\Sigma )$, it follows by the definition
of $L$ that $D(L)\subset H^{1-\varepsilon }_0(\Sigma _+) $. Since $L$
is surjective onto $H^\frac12(\Sigma _+^\circ)$, the last statement follows.
\qed
\enddemo

There is a simple example mentioned in  \cite{S68} of a
harmonic function 
$u(x_1,x_2)=\operatorname{Im}(x_2+ix_1)^{\frac12}$ on ${\Bbb R}_+\times{\Bbb R}$
satisfying the mixed condition on $\{x_1=0\}$, namely $\gamma _0u=0$
for $x_2>0$, $\gamma _1u=0$
for $x_2<0$. It is not in $H^{\frac32}$ in a neighborhood of 0.
This shows that $D(A_{\chi ,\Sigma
_+})$ is not in general contained in $H^{\frac32
}(\Omega )$, so the regularity cannot be improved.

Now consider the Kre\u\i{}n resolvent formula (2.16) for this choice
of $L$ and $X$; by the selfadjointness, $Y=X$.
Recall that  $K _{\gamma ,X}=\inj_{V } \gamma _{V 
}^{-1}\colon X\to L_2(\Omega )$, where $V  $ is the subspace of
$Z  =\ker(\Ama  )$ that is mapped to $X$ by $\gamma _0$. Since $\gamma _{V  }^{-1}$ acts like $K 
_\gamma $ from the space $X= H^{-\frac12}_0(\Sigma _+)$ to $V 
$, we can also write
$$
K _{\gamma ,X}=\inj_{V } K _\gamma e^+
\colon H^{-\frac12}_0(\Sigma _+)\to L_2(\Omega ),\text{ and then } K _{\gamma ,X}^*= r^+
K _\gamma^*\pr_{V  }\colon  L_2(\Omega )\to
H^\frac12(\Sigma _+^\circ),
$$
whereby the formula takes the form
$$
A_{\chi ,\Sigma _+} ^{-1}-A_\gamma  ^{-1}
=\inj_{V  }K _{\gamma }e^+L ^{-1}r^+
K_{\gamma}^*\pr_{V  }
=\inj_{V  }K _{\gamma }e^+(r^+(b-P_{\gamma ,\nu })e^+) ^{-1}r^+
K_{\gamma}^*\pr_{V  }
.\tag4.5
$$

The $\lambda $-dependent version is  found by replacing $A$ by
$A-\lambda $ in the various defining formulas, as explained at the end
of Section 2. Since $\chi =\nu -b\gamma _0$, we have that
$$
 P^\lambda _{\gamma ,\chi }=\chi K^\lambda _\gamma =P^\lambda _{\gamma
,\nu }-b,\tag4.6
$$
a notation we shall now use. Moreover, using the standard
abbreviation for a truncated operator $r^+Qe^+=Q_+$, we can write $r^+P^\lambda _{\gamma ,\chi
}e^+=P^\lambda _{\gamma ,\chi ,+}$. Then the result in the $\lambda
 $-dependent formulation is:

\proclaim{Theorem 4.2} Let $\Sigma _+$ and $b$ be smooth. Then
$$
\aligned
L^\lambda \varphi &=-P^\lambda _{\gamma ,\chi ,+}\varphi \text{ for
}\varphi \in D(L),\, \lambda \in \varrho
(A_\gamma ),\\
(A_{\chi ,\Sigma _+}-\lambda ) ^{-1}-(A_\gamma -\lambda ) ^{-1}
&=-K^\lambda _{\gamma ,X}(P^\lambda _{\gamma ,\chi ,+}) ^{-1}
(K^{\bar\lambda }_{\gamma ,X})^*\text{ for }\lambda \in \varrho (A_{\chi ,\Sigma _+})\cap \varrho
(A_\gamma )
.
\endaligned\tag4.7
$$
where $V_\lambda =K^\lambda _\gamma (X)$, $K^\lambda _{\gamma
,X}=\inj_{V_\lambda }\gamma _{V_\lambda 
}^{-1}=\inj_{V _\lambda  }K ^\lambda _{\gamma }e^+$, $(K^{\bar\lambda
}_{\gamma ,X})^*=(\gamma _{V_{\bar\lambda
}}^{-1})^*\pr_{V_{\bar\lambda }}=r^+(K^{\bar\lambda
}_{\gamma})^*\pr_{V _{\bar\lambda } }$. 
\endproclaim

The inverse of $P^\lambda _{\gamma ,\chi }$ is
$P^\lambda _{\chi ,\gamma }$, when it exists. It is important to
observe that $(P^\lambda _{\gamma ,\chi ,+})^{-1}$ is {\it not} the
same as
$P^\lambda _{\chi ,\gamma ,+}$; this is part of the difficulty treated
in Section 5.

\subhead 4.2 Other Kre\u\i{}n resolvent formulas \endsubhead

Next, if we work instead with a Neumann realization as the 
reference operator, we can show
a different formula containing full Poisson operators.

Consider again the boundary condition
$$
\gamma _0u=0\text{ on }\Sigma _-, \quad \nu u=b\gamma _0u\text{
on }\Sigma _+ .\tag4.8
$$
If $b$ has a bounded inverse $f $, we can set $f
_+=1_{\Sigma _+}f $ and  write condition (4.8) as
one equation, a Dirichlet-type condition
$$
\gamma _0u=f _+\nu u.\tag4.9
$$
Here $\gamma _0u$ is a function of $\nu u$, so that the operator
$A_\nu $ can be used in a simple way as the reference operator.

Actually, it only takes a small modification to obtain invertibility
of the coefficient in general: If $b$ does not have a bounded inverse,
we can replace $\nu u$ by $$
\nu 'u=\nu u+K\gamma _0u,\tag4.10
$$ where
$K$ is chosen $>\operatorname{ess} \sup|b(x)|$ (as in Section 2.1); then the condition
(4.8)
takes the form
$$
\gamma _0u=0\text{ on }\Sigma _-, \quad \nu 'u=b'\gamma _0u\text{
on }\Sigma _+,\tag4.11
$$
where $b'=b+K$ does have a bounded inverse. Note that $\chi '=\nu '-b'\gamma
_0=\chi $ by cancellation. In Green's formula (2.14) we
get $\nu $ replaced by $\nu '$ by adding the term $(K\gamma
_0u,\gamma _0u)-(\gamma _0u,K\gamma _0u)$ (equal to 0) to the
right-hand side, and the
sesquilinear form is adapted to these formulas by addition of 
the first-order terms
$\sum_{j=1}^n[(Kn_j\partial_ju,v)_\Omega +(u,K\partial_j(n_jv))_\Omega ]$, giving the
form
$$
a'(u,v)=a(u,v)+\sum_{j=1}^n[(Kn_j\partial_ju,v)+(u,K\partial_j(n_jv))]\text{
on }H^1(\Omega ).
$$
Here the $n_j$ are extended smoothly to
the interior of $\Omega $, vanishing outside a small neighborhood of
$\Sigma $. The operators defined from $a'$ on various spaces between $H^1(\Omega
)$ and $H^1_0(\Omega )$ still act like $A$, since
$(u,K\partial_j(n_j\varphi ))=-(Kn_j\partial_ju,\varphi )$ for
$\varphi \in C_0^\infty (\Omega )$. The ``halfways Green's formula''
is here
$$
(Au,v)-a'(u,v)=(\nu 'u,\gamma _0v)_{L_2(\Sigma )} ,\tag4.12
$$
since $
\sum_{j=1}^n[(Kn_j\partial_ju,v)_{\Omega }+(u,K\partial_j(n_jv ))_{\Omega
}]=-(K\sum n_j^2\gamma _0u,\gamma _0v)_{\Sigma }=-(K\gamma _0u,\gamma _0v)_\Sigma .$
The forms in the scheme (2.2) are now replaced by
$$
\aligned
a'_\gamma (u,v)&=a'(u,v)\text{ on }H^1_0(\Omega ) ,\text{ leading to }A_\gamma ,\\
a'_{\nu '} (u,v)&=a'(u,v)\text{ on }D(a_{\nu '} )=H^1(\Omega ),\text{
leading to }A_{\nu '},\\
a'_{\chi '} (u,v)&=a'(u,v)+(b'\gamma _0u,\gamma _0v)_{L_2(\Sigma )}\text{ on }D(a_{\chi '}  )=H^1(\Omega ),\text{
leading to }A_{\chi } ,\\
a'_{\chi ',\Sigma _+} (u,v)&=a'(u,v)+(b'\gamma _0u,\gamma _0v)_{L_2(\Sigma
_+)}\text{ on }
 D(a'_{\chi ',\Sigma _+}  )=H^1_{\Sigma _+}(\Omega )\\& =\{u\in H^1(\Omega )|\supp u\subset
\Sigma _+\},\text{
leading to }A_{\chi ,\Sigma _+} ;
\endaligned\tag 4.13
$$
here $A_{\nu '}$ is the realization of $A$ under the boundary
condition $\nu 'u=0$, whereas the choices with $b'$ still give the
boundary condition $\nu u=b\gamma _0u$ on $\Sigma $ resp.\ $\Sigma
_+$, since $b'=b+K$, $\nu '=\nu +K\gamma _0$. With $K _{\nu
'}$, $P_{\gamma ,\nu ' } $ and $P_{\nu ',\gamma  } $
defined as in Section 4.1 with $\nu $ replaced by
$\nu '$, and 
$$
\Gamma '_{\gamma }= \gamma _0 -P_{\nu ',\gamma  }\nu ',
$$
 we have the generalized Green's formula valid for $u,v\in D(\Ama)$:
$$
(Au,v)_{L_2(\Omega )}-(u,Av)_{L_2(\Omega
)}=
(\nu 'u,\Gamma ' _{\gamma }v)_{-\frac32,\frac32}-(\Gamma
'_{\gamma } u,\nu 'v)_{\frac32,-\frac32}.
\tag4.14$$

In the following, we assume that the forms in (4.13) all have positive
lower bound. We set $f =(b')^{-1}$ so that the mixed boundary
condition (4.11) can be written
$$
\gamma _0u=f _+\nu 'u,\tag4.15
$$
where $f _+= 1_{\Sigma _+}f $, as accounted for above.

We now describe $A_{\chi ,\Sigma _+}$ in terms of the correspondence (2.12) and
its interpretation in Section 2.4, with $\nu $ replaced by $\nu '$.

Here $X_1$ is
the full space $H^{-\frac32}(\Sigma )$, which is seen as follows:
When $\psi \in C_0^\infty (\Sigma _+^\circ)\cup C_0^\infty (\Sigma
_-^\circ)$, then $f _+ \psi \in C_0^\infty (\Sigma
_+^\circ)$, and there exists $u\in C^\infty (\comega)$ such that ${\nu '}
u=\psi $, $\gamma _0u = f _+ \psi $; this $u$ satisfies
(4.9). So $ C_0^\infty (\Sigma _+^\circ)\cup C_0^\infty (\Sigma
_-^\circ)\subset D(L _1)$. It is known that $ C_0^\infty (\Sigma _+^\circ)\cup C_0^\infty (\Sigma
_-^\circ)$ is dense in $H^{s}(\Sigma )$ for $s<\frac12$.
In particular $ C_0^\infty (\Sigma _+^\circ)\cup C_0^\infty (\Sigma
_-^\circ)$ is dense in $H^{-\frac32}(\Sigma )$, so we conclude that
$X_1=H^{-\frac32}(\Sigma )$. Since $A_{\chi ,\Sigma _+}$ is selfadjoint, also
$Y_1=H^{-\frac32}(\Sigma )$.

Thus the realization $A_{\chi ,\Sigma _+}  $ with domain (2.6) corresponds 
to an operator $L_1  
\colon H^{-\frac32}(\Sigma )\to H^{\frac32}(\Sigma )$ with domain
$D(L_1  )={\nu '} D(A_{\chi ,\Sigma _+})$; the latter lies in
$H^{-\frac12}(\Sigma )$ since $D(A_{\chi ,\Sigma _+})\subset H^1(\Omega )$.
It follows by comparison of (2.30) with (4.15) that $L_1$ acts as
$$
L _1 =-f _+ + P _{{\nu '} ,\gamma }. \tag4.16
$$
Since $A_{\chi ,\Sigma _+}$ is bijective, so is $L_1$.

Then the Kre\u\i{}n resolvent formula reads
$$
A_{\chi ,\Sigma _+} ^{-1}-A_{\nu '}   ^{-1}=K _{{\nu '} }L_1  ^{-1}K_{\nu '} ^*=K
_{{\nu '} }
(P _{{\nu '} ,\gamma }-f _+ )^{-1}K_{\nu '} ^*,\text{ for }\lambda \in
\varrho (A_{\chi ,\Sigma _+})\cap \varrho (A_{\nu '}).\tag4.17
$$
It may look a little more useful than (4.5), since the operators
surrounding $L _1^{-1}$ are a full Poisson operator
and trace operator in the pseudodifferential boundary operator
calculus, but it poses again the question of a detailed understanding of the
term in the middle, defined on $H^\frac32(\Sigma )$. This may not
be any easier than our treatment in Section 4.1, since the principal
part of $L_1$ is the 0-order multiplication by $-f _+ $ which
vanishes on $\Sigma _-$, and $P _{{\nu '}
,\gamma }$ is of order $-1$.

We can replace $A$ by $A-\lambda $ in the various defining formulas
and obtain:

\proclaim{Theorem 4.3} Let $\Sigma _+$ and $b$ be smooth. Define $\nu '$ by {\rm (4.10)ff.}\ and $f
=(b+K)^{-1}$. Then 
$$
(A_{\chi ,\Sigma _+}-\lambda ) ^{-1}-(A_{\nu '}  -\lambda ) ^{-1}=K^\lambda
_{{\nu '} }(P^\lambda  _{{\nu '} ,\gamma }-f _+
)^{-1}(K^{\bar\lambda }_{\nu '}) ^*,\text{ for }\lambda \in \varrho
(A_{\nu '})\cap \varrho (A_{\chi ,\Sigma _+}).\tag4.18
$$
\endproclaim

Formula (4.18) can even be turned into a resolvent difference formula where the surrounding
Poisson operator and trace operator come from the Dirichlet problem,
by use of the fact that $P ^\lambda _{\gamma ,{\nu '} }$ and $P^\lambda  
_{{\nu '} ,\gamma }$ are inverses of one another, and
$$
K ^\lambda _{\nu '} =K ^\lambda _\gamma  P ^\lambda _{{\nu '} ,\gamma  },\quad (P^{\bar\lambda } _{{\nu '} ,\gamma })^*=P ^\lambda _{{\nu '} ,\gamma },\tag4.19
$$
then. Namely, insertion in (4.18) gives:
$$
(A_{\chi ,\Sigma _+}-\lambda ) ^{-1}-(A_{\nu '} -\lambda )  ^{-1}=K ^\lambda 
_{\gamma   }P^\lambda  _{{\nu '} ,\gamma  }(P^\lambda  _{{\nu '} ,\gamma }-f _+ )^{-1}P^\lambda _{{\nu '} ,\gamma }(K^{\bar\lambda } _{\gamma  })^*.\tag4.20
$$
This can be added to the well-known  formula 
$$
(A_{{\nu '} }-\lambda ) ^{-1}-(A_\gamma -\lambda ) ^{-1}=K^\lambda  _{\gamma }(-P^\lambda  _{\gamma ,\nu '})^{-1}(K^{\bar\lambda } _{\gamma })^*=-K^\lambda  _{\gamma }P^\lambda  _{{\nu '} ,\gamma  } (K^{\bar\lambda } _{\gamma })^*
$$
((2.32) with
$\wA=A_{\nu '}$, hence $L^\lambda =-P^\lambda _{\gamma ,\nu '}$),
to give a formula for the resolvent difference with the
Dirichlet realization, having another structure than (4.7):
$$
\multline
(A_{\chi ,\Sigma _+}-\lambda ) ^{-1}-(A_\gamma  -\lambda )  ^{-1}
=K ^\lambda 
_{\gamma   }P^\lambda  _{{\nu '} ,\gamma }[(P^\lambda  
_{{\nu '} ,\gamma }-f _+ )^{-1}P^{\lambda }_{{\nu '} ,\gamma  }-1](K^{\bar\lambda }_{\gamma  })^*\\
=K^\lambda  
_{\gamma   }P^\lambda  _{{\nu '} ,\gamma }(P^\lambda  
_{{\nu '} ,\gamma }-f _+ )^{-1}f _+ (K^{\bar\lambda }_{\gamma  })^*
, \text{ for }\lambda \in
\varrho (A_{\chi ,\Sigma _+})\cap \varrho (A_\gamma ) \cap \varrho (A_{\nu '} ) .
\endmultline\tag4.21
$$

The last formula in (4.21) has a similar flavor as the formula found
by Pankrashkin in \cite{P06}, Sect.\ 4.3.
 
If $b$ itself is invertible, the formulas will be valid with $f
=b^{-1}$, $\nu '$ replaced by $\nu $. We have shown:

\proclaim{Corollary 4.4} Under the hypotheses of Theorem {\rm 4.3}, we
have the formulas in {\rm (4.21)} for the difference with the
Dirichlet resolvent, when $\lambda \in
\varrho (A_{\chi ,\Sigma _+})\cap \varrho (A_\gamma ) \cap \varrho (A_{\nu '} )$.

If $b$ itself is invertible, there are the formulas with $f
=b^{-1}$:
$$\align
(A_{\chi ,\Sigma _+}-\lambda ) ^{-1}-(A_{\nu }  -\lambda ) ^{-1}&=K^\lambda
_{{\nu } }(P^\lambda  _{{\nu } ,\gamma }-f _+
)^{-1}(K^{\bar\lambda }_{\nu }) ^*,\tag4.22\\
(A_{\chi ,\Sigma _+}-\lambda ) ^{-1}-(A_\gamma  -\lambda )  ^{-1}
&=K^\lambda  
_{\gamma   }P^\lambda  _{{\nu } ,\gamma }(P^\lambda  
_{{\nu } ,\gamma }-f _+ )^{-1}f _+ (K^{\bar\lambda }_{\gamma  })^*
, \tag4.23
\endalign$$
where {\rm (4.22)} holds  for $\lambda \in \varrho
(A_{\nu })\cap \varrho (A_{b ,\Sigma _+})$, {\rm (4.23)} holds for $\lambda \in
\varrho (A_{\chi ,\Sigma _+})\cap \varrho (A_\gamma ) \cap \varrho (A_{\nu } )$.
\endproclaim

\example{Remark 4.5} The analysis in Proposition 4.1ff.\ showed that
$D(A_{\chi ,\Sigma _+})\subset H^{\frac32-\varepsilon }(\Omega )$ but is
not in general
contained in $H^\frac32(\Omega )$. Thus those results in Malamud
\cite{M10}, Section 6, that concern second-order realizations with domain contained  in
$H^\frac32(\Omega )$ (i.e., with $\gamma _0u$ and $\nu u\in L_2(\Sigma )$), 
will not in general apply
to the mixed problem.
\endexample

\example{Remark 4.6}
If we instead of (4.9) consider a boundary condition 
$$
\gamma _0u=g {\nu } u,\tag4.24
$$
 where $g $ is an arbitrary  $C^\infty $-function on $\Sigma $,
we can carry an analysis through, showing that if the corresponding
realization $\wA$ is bijective and selfadjoint, then it corresponds to
an operator $L_1$ from $H^{-\frac32}(\Sigma )$ to $H^\frac32(\Sigma
)$, with domain dense in $H^{-\frac32}(\Sigma )$ and acting like $P_{{\nu }
,\gamma }-g $, such that there are Kre\u\i{}n formulas$$\aligned
\wA ^{-1}-A_{\nu }   ^{-1}&=K _{{\nu } }L_1  ^{-1}{K_{\nu }} ^*=K
_{{\nu } }(P _{{\nu } ,\gamma }-g  )^{-1}{K_{\nu }} ^*,\\
\wA ^{-1}-A_\gamma    ^{-1}
&=K 
_{\gamma   }P _{{\nu } ,\gamma }
(P _{{\nu } ,\gamma }-g  )^{-1}g  {K_{\gamma  }}^*,
\endaligned\tag4.25$$
and $\lambda $-dependent variants. But again, the operator $L_1^{-1}=(P _{{\nu } ,\gamma }-g  )^{-1}$ is nonstandard in
the calculus of $\psi $do's, since  $P _{{\nu } ,\gamma }$ is elliptic
of order $-1$ whereas $g $ defines an operator of order 0 and can
vanish on large subsets of $\Sigma $.
\endexample

\subhead 4.3 Nonsmooth domains 
\endsubhead

We here include some observations on cases where the set $\Omega $ is
not smooth. An interesting variant of the Zaremba problem is
where $\Sigma =\Sigma _+\cup\Sigma _-$ with $\Sigma _+$ and $\Sigma _-$
{\it meeting at  an angle $< \pi $}. Then there is the perhaps surprising fact
that the realization $A_{\nu ,\Sigma _+}$ of $-\Delta $ with Dirichlet condition $\gamma _0u=0$
on $\Sigma _-$ and Neumann condition $\nu u=0$ on $\Sigma _+$ ($\nu
=\gamma _1$) can 
have a {\it better} regularity than when $\Omega $ is
smooth. Here is an example:

\example{Example 4.7} Let $ \Omega '$ be a smooth bounded set
that is symmetric in $x_1$ around $x_1=0$; i.e., is preserved under the
mapping $J_1\colon (x_1,x_2,\dots, x_n)\to (-x_1,x_2,\dots, x_n)$. Let $\Omega
=\{x\in \Omega ' \mid x_1>0\}$. Then the solutions of the mixed problem
for $-\Delta $ on $\Omega $ with $\Sigma _-=\partial\Omega '\cap \{x_1\ge 0\}$,
$\Sigma _+=\overline{\Omega'}\cap \{x_1=0\}$, are the restrictions to
$\Omega $ of those solutions to the Dirichlet problem for $\Omega '$ 
that are invariant
under $J_1$. (This observation enters in a prominent way in the
discussion of isospectral domains for mixed problems by Levitin,
Parnovski and Polterovich \cite{LPP06}.) Here the domain of the
Dirichlet realization of $-\Delta $ on $\Omega '$ is in $H^2(\Omega
')$, hence $D(A_{\nu ,\Sigma _+})\subset H^2(\Omega )$ (observe that both
operators are bijective when defined by the variational construction).
In this case $\Sigma _+$ and $\Sigma _-$ meet at an angle $\pi/2$.
--- Related results are found for polygonal domains, cf.\ Grisvard \cite{G85}.
\endexample

More generally, consider the case where
$\Omega $ is such that 
$ \Sigma _+$ and $\Sigma _-$ meet at an angle $
<\pi $, in the way described in Brown \cite{B94}; such domains are
by some authors called {\it creased domains}. It is shown there that the solutions $u\in H^1(\Omega )$ of 
$$
-\Delta u=0\text{ in }\Omega ,\quad \gamma _0u=\varphi \text{ on }\Sigma
_-,\quad \gamma _1u=\psi  \text{ on }\Sigma _+,\tag4.26
$$
with $\varphi \in H^1(\Sigma ^\circ _-)$, $\psi \in L_2(\Sigma _+)$, have $\gamma
_0(\nabla u)\in L_2(\Sigma )$; in particular $\gamma _1u\in L_2(\Sigma
)$. Here $\Sigma $ just needs to be Lipschitz, in such a way that
$\partial\Sigma _+$ is Lipschitz in $\Sigma $ (we refer to \cite{B94}
for the precise description).           

To apply this to $A_{\nu ,\Sigma _+}$, we
restrict to quasi-convex domains $\Omega $. They are defined by 
Gesztesy and Mitrea in \cite{GM11} as
a special case of Lipschitz domains including convex domains, which
allow showing solvability and regularity 
theorems for the Dirichlet and Neumann problems for
$-\Delta $ on $\Omega $ in larger  scales of Sobolev-type spaces than
in Jerison and Kenig \cite{JK95}; the work builds on Mitrea, Taylor and Vasy
\cite{MTV05} and Mazya, Mitrea and Shaposhnikova \cite{MMS10}.

\proclaim{Theorem 4.8} Assume that $\Omega $ is bounded, open and
quasi-convex as defined in {\rm \cite{GM11}}. Assume moreover that
$\Omega $ is creased, in the way that the boundary
$\Sigma $ equals $\Sigma _+\cup\Sigma _-$, where $\Sigma _+$ and $\Sigma _-$
meet at an angle $<\pi $, as described in
{\rm \cite{B94}}. The realization $A_{\nu ,\Sigma _+}$ of $-\Delta $ with
Neumann condition on $\Sigma _+$, Dirichlet condition on $\Sigma _-$
then has $D(A_{\nu ,\Sigma _+})\subset H^{\frac32}(\Omega )$. 
\endproclaim

\demo{Proof} We can assume that $\Sigma _-\ne \emptyset$.
To describe a solution in $H^1(\Omega )$ of 
$$
-\Delta u=f\text{ in }\Omega ,\quad \gamma _0u=0\text{ on }\Sigma
_-,\quad \gamma _1u=0 \text{ on }\Sigma _+,\tag4.27
$$
with $f\in L_2(\Omega )$, let $v$ be the solution of the Dirichlet
problem
$$
-\Delta v=f\text{ in }\Omega ,\quad \gamma _0v=0\text{ on }\Sigma
;\tag4.28
$$
then $z=u-v$ should be a solution of (4.26) with $\varphi =0$, $\psi
=-\gamma _1v|_{\Sigma _+}$.
Since $v\in H^2(\Omega )$ by \cite{GM11} Th.\ 10.4, $\gamma _1v|_{\Sigma _+}\in
H^\frac12(\Sigma ^\circ _+)\subset L_2(\Sigma _+)$, so the result of
Brown \cite{B94} implies that $\gamma _1z\in H^1(\Sigma )$. Then the
regularity theorem for the Neumann problem \cite{GM11} Th.\ 10.8
implies that $z\in H^\frac32(\Omega )$, hence $u=v+z\in H^\frac32(\Omega
)$. (Since $A_{\nu ,\Sigma _+}$ is bijective, the solutions we
consider are consistent with those considered in \cite{B94}.)  We conclude that
$D(A_{\nu ,\Sigma _+})\subset H^\frac32(\Omega )$.\qed 
\enddemo

Note the contrast with the informations
obtained in Proposition 4.1ff.\,  where $\Sigma $ is smooth and $D(A_{\nu ,\Sigma _+})$ is
in general only in $H^{\frac32-\varepsilon }(\Omega )$. But that is a case
where $\Sigma _+$ and $\Sigma _-$ meet at the angle $\pi $, which is
explicitly excluded in \cite{B94}.

The mixed problem in these various forms can, for piecewise smooth domains,
be regarded as a special
case of crack problems and edge problems, as studied e.g.\ by
Duduchava, Dauge, Costabel, Mazya, Solonnikov and their collaborators, 
see also Schulze et
al.\ \cite{RS83, HS08}. The results are often described in terms of
norms weighted by powers of the distance to the edge; this gives a
clarification of
the singularities, but can lead outside the Sobolev spaces considered here.

Let us finally mention that for quasi-convex domains 
there is in \cite{GM11} Th.\ 10.4 
established a homeomorphism $\hat\gamma _D\colon
Z\simto (N^{\frac12}(\Sigma ))^*$ (generalizing $\gamma _0$), which allows translating formula
(2.10) in Section 2 above to a formula like (2.16) with (2.17). 
Here $N^{\frac12}(\Sigma )$ is a certain Hilbert space related to
$H^\frac12(\Sigma )$ explained in \cite{GM11}, and
$(N^{\frac12}(\Sigma ))^*$ is its dual space with respect to a
sesquilinear duality consistent with the $L_2(\Sigma )$-scalar
product, such that
$
N^{\frac12}(\Sigma )\subset L_2(\Sigma )\subset (N^{\frac12}(\Sigma ))^*,
$
with dense, continuous injections. For a general closed realization
$\wA$ of $-\Delta $, let $X$ and $Y$ be the closures of
$\hat\gamma _D(D(\wA))$ resp.\ $\hat\gamma _D(D(\wA^*))$ 
in $(N^{\frac12}(\Sigma ))^*$, and let $V$ resp.\ $W$ be their
inverse images in $Z$ (by  $\hat\gamma _D^{-1}$); in the
selfadjoint case, $Y=X$.
We can then define $\gamma _V$ to be the
restriction of $\hat\gamma _D$ mapping $V$ homeomorphically to
$X$; similarly, $\gamma _W$ is the
restriction of $\hat\gamma _D$ mapping $W$ homeomorphically to
$Y$. 
With this, the considerations in (2.15)--(2.18) are valid, leading to:

\proclaim{Theorem 4.9} When $\Omega $ is 
quasi-convex and $\wA$ is a general closed realization of $-\Delta $
with $0\in \varrho (\wA)$,
it satisfies the Kre\u\i{}n resolvent formula {\rm (2.16)} with {\rm
(2.17)}, $L$ defined as after {\rm (2.15)}.
In particular, for the realization $A_{\nu ,\Sigma
_+}$ of the mixed problem,  one has with $X=$ the closure of $\hat\gamma _D(D(A_{\nu ,\Sigma _+}))$,
$$A_{\nu ,\Sigma _+}^{-1}
-A_\gamma ^{-1}=K _{\gamma ,X} 
L^{-1}(K_{\gamma ,X})^*,\quad K _{\gamma
  ,X}= \inj_{V  }\gamma _{V }^{-1}\colon X\to
V \subset H.\tag4.29
$$
\endproclaim 

The rest of Section 2 likewise carries over to the quasi-convex setting,
but it must be noted that the Dirichlet-to-Neumann operator is then an
abstractly 
defined operator whose local structure is not so well known. Also
$\lambda $-dependent variants of (2.16) for realizations of $-\Delta
-\lambda $ are valid when $\lambda \in \varrho (A_\gamma )\cap \varrho
(\wA)$. The interpretation of the general theory of
\cite{G68} for quasi-convex domains is 
worked out in great detail in \cite{GM11}. 

We remark however that the Kre\u\i{}n formula in \cite{GM11}
Th.\ 16.3  differs from our formula (2.16), particularly when $X\ne
(N^{\frac12}(\Sigma ))^*$ (which holds for  genuine mixed problems).

Upper eigenvalue estimates (1.6)
for the resolvent difference (4.29) follow from \cite{B62}, 
cf.\ Remark 3.3. Asymptotic estimates would demand an effort that to our
knowledge has not yet been taken up.

\head 5. Spectral asymptotics for the mixed problem\endhead

\subhead 5.1 Notation \endsubhead

In this section we go back to smooth domains and restrict the attention to the case $a_{jk}=\delta
_{jk}$, i.e., we take  $A$  principally equal to $-\Delta $, 
in order to use some detailed formulas in Eskin \cite{E81}.

We want to show a spectral asymptotic
formula for the operator
$$
(A_{\chi ,\Sigma _+}-\lambda ) ^{-1}-(A_\gamma -\lambda ) ^{-1}
=K^\lambda _{\gamma ,X}(L^\lambda ) ^{-1}
(K^{\bar\lambda }_{\gamma ,X})^*=-K^\lambda _{\gamma ,X}(P^\lambda _{\gamma ,\chi ,+}) ^{-1}
(K^{\bar\lambda }_{\gamma ,X})^*
$$ from Theorem 4.2. As done also earlier, we begin by taking 
$\lambda $ as a sufficiently low fixed real 
number such that the considered realizations of $A-\lambda $ are
positive, and then omit $\lambda $ from the notation. General $\lambda
$ are included in the proof of the final Theorem 5.17.

In view of the
formula (2.17) for $K_{\gamma ,X }$,
we are considering the operator
$$
A_{\chi ,\Sigma _+} ^{-1}-A_\gamma  ^{-1}=\inj_V\gamma _V^{-1}L^{-1}(\gamma _V^{-1})^*\pr_V;\tag5.1
$$
it is compact selfadjoint nonnegative.

Let us first recall some facts on spaces describing the spectral behavior
of compact operators.
For $p>0$ we denote by $\Cal C_p$ the Schatten class of compact linear
operators $B$ (in a Hilbert space $H$) with singular value sequences
$(s_j(B))_{j\in{\Bbb N}}$ belonging to $
\ell_p$, and by $\frak S_{p}$ the quasi-normed space of
compact operators $B$ with $s_j(B) =O(j^{-1/p})$ (sometimes called a
weak Schatten class); here $\frak
S_{p}\subset \Cal C_{p' }$ for $p'>p$. 
Moreover, we denote by $\frak S_{p,0}$ the subset of operators $B\in
\frak S_p$ for which $s_j(B) =o(j^{-1/p})$, i.e., $s_j(B)j^{1/p}\to 0$
for $j\to\infty $. Clearly, $\frak S_{p}\subset\frak S_{p',0}$ for $p'>p$.

The rules shown by Ky Fan \cite{F51}$$
s_{j+k-1}(B+B')\le s_j(B)+s_k(B'),\quad s_{j+k-1}(BB')\le
s_j(B)s_k(B'),
$$ 
imply that $\Cal C_p$, $\frak S_p$ and  $\frak S_{p,0}$ are vector spaces, and that
there are the following
product rules:$$
 \Cal
C_p\cdot \Cal C_q\subset \Cal C_{1/(p^{-1}+q^{-1})}, \quad \frak
S_p\cdot \frak S_q\subset \frak S_{1/(p^{-1}+q^{-1})},\quad \frak S_p\cdot \frak S_{q,0}\subset \frak S_{1/(p^{-1}+q^{-1}),0}.\tag5.2
$$
Moreover, the rule for $F_1,F_2\in \Cal L(H)$, $$
s_j(F_1BF_2)\le \|F_1\|s_j(B)\|F_2\|\tag5.3
$$ implies that $\Cal C_p$, $\frak
S_p$ and $\frak S_{p,0}$ are preserved under compositions with bounded
operators. They are also preserved under taking adjoints. We recall two perturbation results:

\proclaim{Lemma 5.1}

$1^\circ$ If $s_j(B)j^{1/p}\to C_0$ and $s_j(B')j^{1/p}\to 0$ for
$j\to\infty$, then $s_j(B+B')j^{1/p}\to C_0$ for $j\to\infty$.

$2^\circ$ If $B=B_M+B'_M$ for each $M\in\Bbb N$, where
$s_j(B_M)j^{1/p}\to C_M$ for $j\to\infty$ and $s_j(B'_M)j^{1/p}\le
c_M$ for
$j\in\Bbb N$, with $C_M\to C_0$ and $c_M\to 0$ for
$M\to\infty$, then $s_j(B)j^{1/p}\to C_0$ for $j\to\infty$.
\endproclaim

The statement in $1^\circ$ is the Weyl-Ky Fan theorem (cf.\ e.g.\
\cite{GK69} Th.\ II 2.3), and $2^\circ$ is a refinement shown in
\cite{G84}, Lemma 4.2.$2^\circ$.

We also recall that when $\Xi$ and $\Xi _1$ are $m$-dimensional 
manifolds (possibly with a boundary, sufficiently smooth), $\overline \Xi _1$ being compact,  and $B$ is a bounded linear operator
from $L_2(\Xi)$ to $H^t(\Xi_1)$ for some $t>0$, then $B\in \frak
S_{m/t}$ as an operator from $L_2(\Xi)$ to $L_2(\Xi _1)$, with $$
s_j(B)j^{t/m}\le C\|B\|_{\Cal L(L_2(\Xi ),H^t(\Xi _1))},\tag 5.4
$$
with a constant $C$ depending on $t$ and the manifolds (references
e.g.\ in \cite{G84}).

\subhead 5.2 Constant coefficients \endsubhead

One ingredient in the analysis of the spectrum of (5.1) is an application of the
constant-coefficient situation, so we begin by working that out, in the
case
$b=0$.
Here
$\Omega $,  $\Sigma $ and $\Sigma _\pm$ are replaced by $\rnp$,  ${\Bbb
R}^{n-1}$ and $\crpm^{n-1}$, and we take $A=-\Delta +\alpha ^2$ for
some $\alpha >0$; marking the operators with a subscript 0.
The Poisson operator $K_{0,\gamma }$ solving the Dirichlet problem is
the operator $\varphi (x') \mapsto \Cal F^{-1}_{\xi '\to
x'}[e^{-x_n(|\xi '|^2+\alpha ^2)^{\frac12}}\hat\varphi (\xi ')]$, so
the Dirichlet-to-Neumann operator $P_{0,\gamma ,\nu }$ is the $\psi
$do with symbol $-(|\xi '|^2+\alpha ^2)^{\frac12}$, i.e.,
$$
P_{0,\gamma ,\nu }=-\Op((|\xi '|^2+\alpha ^2)^\frac12)=-(-\Delta
_{x'}+\alpha  ^2)^{\frac12}, \text{ with inverse }P_{0,\nu ,\gamma }=-\Op((|\xi '|^2+\alpha ^2)^{-\frac12}).
$$
($\Cal F$ denotes the Fourier transform, and $\Op(a(x',\xi '))v=\Cal
F^{-1}_{\xi '\to x'}(a(x',\xi ')\Cal F_{x'\to\xi '}v)$.) Then with $\xi ''=(\xi _1,\dots,\xi
_{n-2})$,
$$
L_0  =-r^+P_{0,\gamma ,\nu }e^+=r^+\Op((|\xi ''|^2+\xi _{n-1}^2+\alpha ^2
)^{\frac12})e^+
\colon H^s_0(\rp^{n-1})\to H^{s-1}(\rp^{n-1});
%, \quad 0<s<1
\tag5.5$$
it will be used with $s=1-\varepsilon $, cf.\ Proposition 4.1.
According to Eskin \cite{E81}, Ch.\ 7, one has in 
view of the factorization
$$
(|\xi ''|^2+\xi _{n-1}^2+\alpha ^2
)^{\frac12}=((|\xi ''|^2+\alpha ^2 )^\frac12-
i\xi_{n-1})^{\frac12}((|\xi ''|^2+\alpha ^2 )^\frac12 +
i\xi_{n-1})^{\frac12},
$$
that $L_0 $ has the inverse
$$
L_0  ^{-1}= r^+\Lambda _+e^+r^+\Lambda _-\ell^+_s
\colon H^{s-1}(\rp^{n-1})\to
H^s_0(\rp^{n-1}),\quad 0<s<1, \tag5.6$$
where 
$$\Lambda _\pm=\Op(\lambda _\pm(\xi ')),\quad 
\lambda _\pm(\xi ')=((|\xi ''|^2+\alpha ^2 )^\frac12\pm
i\xi_{n-1})^{-\frac12},
\tag5.7
$$ 
and $\ell^+_s$ denotes a smooth extension operator, continuous from
$H^{t}(\rp^{n-1})$ to $H^t({\Bbb R}^{n-1} )$ for all $t$.
The operators $\Lambda _\pm$
are a ``plus-operator'' resp.\ a
``minus-operator'' in the terminology of \cite{E81}; plus-operators
preserve support in $\crp^{n-1}$, and minus-operators are adjoints of
plus-operators and preserve support in $\crm^{n-1}$. 

When the formula is used for $s=1-\varepsilon  >\frac12$, we
can replace $\ell^+_s$ by $e^+$, so
$$
L_0 ^{-1}= r^+\Lambda _+e^+r^+\Lambda _-e^+=\Lambda _{+,+}\Lambda _{-,+}\colon H^{-\varepsilon }(\rp^{n-1})\to
H^{1-\varepsilon }_0(\rp^{n-1})\tag 5.8
$$
(recall the notation $Q_+=r^+Qe^+$). $L_0^{-1}$ is of course
different from $(\Lambda _+\Lambda 
_-)_+=-P _{0,\nu ,\gamma ,+}$, that we shall compare it
with further below. We note that $\Lambda _{-,+}$ maps $H^{-\varepsilon
}(\rp^{n-1})$ to $H^{\frac12-\varepsilon }(\rp^{n-1})=
H^{\frac12-\varepsilon }_0(\rp^{n-1})$. Then the fact that $\Lambda _{+,+}$
preserves support in $\crp^{n-1}$, confirms that the range of $L^{-1}$
is in the subspace $H^{1-\varepsilon }_0(\rp^{n-1})$ of $H^{1-\varepsilon }(\rp^{n-1})$.

We shall treat our general problem by reducing to cases in local
coordinates with ingredients principally of this form. 
Then $L_0 ^{-1}$ is
multiplied on both sides with 
cutoff functions, 
so we shall now also consider $\psi L_0 ^{-1}\psi _1$,
where $\psi ,\psi _1\in C_0^\infty (B_R)$ for some ball
$B_R=\{|x'|<R\}\subset {\Bbb R}^{n-1}$. It is continuous 
$$
\psi L_0 ^{-1}\psi _1\colon L_2(B_R\cap \rp^{n-1})\to
H^{1-\varepsilon }(B_R\cap \rp^{n-1}), \text{ any }\varepsilon >0;\tag 5.9
$$
hence in view of (5.4),
$$
\psi L_0 ^{-1}\psi _1\in \frak S_{n-1+\delta },\text{ any }\delta >0.\tag5.10$$ (Better estimates will be
obtained below.) We shall compare it with $-\psi r^+P_{0,\gamma ,\nu
}^{-1}e^+\psi _1=-\psi P_{0,\nu ,\gamma ,+}\psi _1$, and for this purpose we observe that
$$
\aligned
-r^+P_{0,\nu ,\gamma  }e^+-L ^{-1}&= r^+\Lambda _+\Lambda _-e^+- r^+\Lambda _+e^+r^+\Lambda _-e^+\\
&= r^+\Lambda _+e^-JJr^-\Lambda _-e^+=G^+(\Lambda _+)G^-(\Lambda _-),
\endaligned
$$
 where $r^-$ is the restriction operator from ${\Bbb
R}^{n-1}$ to  ${\Bbb
R}^{n-1}_-$,  $e^{-}$ is the corresponding extension-by-zero
operator,
and $J$ is the
reflection operator $J\colon u(x'',x_{n-1})\mapsto u(x'',-x_{n-1})$. We have
used that $I-e^+r^+=e^-r^-$, and denoted
$$
G^+(Q)=r^+Qe^-J,\quad G^-(Q)=Jr^-Qe^+,\tag5.11
$$
as in \cite{G84} and subsequent papers and books of the author. Note that the
distribution kernel of $G^+(\Lambda _+)$ is obtained from that 
of $\Lambda _+$ by
restriction to the second quadrant in $(y_{n-1},x_{n-1})$-space, so
that the singularity at the diagonal $\{x_{n-1}=y_{n-1}\}$ is only felt at 0. 

On the manifold $\Sigma =\Sigma _+\cup\Sigma _-$, $G^{\pm}(Q)$ make
sense only in local coordinates, but
$$
L(Q_1,Q_2)=(Q_1Q_2)_+-Q_{1,+}Q_{2,+}\tag5.12
$$
is well defined when $Q_1$ and $Q_2$ are of order $\le 0$, and locally
has the structure $G^+(Q_1)G^-(Q_2)$.

For later purposes we recall the result of Laptev \cite{L81} (also
shown for $\psi $do's having the transmission property in \cite{G84}):

\proclaim{Theorem 5.2} {\rm \cite{L81}} Let $n-1\ge 2$. When $Q$ is a $\psi $do on ${\Bbb R}^{n-1}$ of
order $-r<0$, and $\psi \in C_0^\infty ({\Bbb R}^{n-1})$, then $\psi G^\pm(Q)$ and $G^{\pm}(Q)\psi $ are in
$\frak S_{(n-2)/r}$, with $s_jj^{r/(n-2)}$ converging to a limit
determined from the principal symbol.  

When $Q_1$ and $Q_2$ are $\psi $do's on $\Sigma =\Sigma _+\cup\Sigma _-$ of
orders $-r_1,-r_2<0$, then $L(Q_1,Q_2)$ is in
$\frak S_{(n-2)/(r_1+r_2)}$.
\endproclaim

The operators $\Lambda _\pm$ are of order $-\frac12$, but are not standard $\psi $do's,
since the symbols $\lambda _\pm$ are not in H\"o{}rmander's symbol
space $S^{-\frac12}_{1,0}$ as
functions of $\xi '$ (high derivatives in $\xi ''$ do
not satisfy the required estimates in terms of powers of $1+|\xi
'|$). Then Laptev's theorem is not applicable to $G^\pm(\Lambda
_+)$ and $G^\pm(\Lambda _-)$. In fact, one
can check that the associated integral operator kernels, calculated
explicitly,  do not satisfy
all the estimates required for Th.\ 3 in \cite{L81}. 
We expect that it should be possible to show  a
spectral estimate as in Theorem 5.2 for these operators, but
leave out further 
investigations here,  settling for some weaker
estimates that still serve our purpose. 

In the following, we denote $x_{n-1}=t$, $y_{n-1}=s$, with dual variables
$\tau ,\sigma $, to simplify the notation. Let $\zeta (t)\in C^\infty ({\Bbb R})$, taking
values in $[0,1]$ and equal to $1$ for $t\ge 1$, equal to $0$ for
$t\le \frac23$. For $\varepsilon >0$, denote $\zeta (t/\varepsilon
)=\zeta _\varepsilon (t)$. 

\proclaim{Lemma 5.3} Let $\varepsilon  >0$. 
The operators $\zeta
_\varepsilon G^+(\Lambda _\pm)$ are of order $-\frac 32$, and $\psi \zeta_\varepsilon
G^+(\Lambda _\pm )$ as well as $\zeta_\varepsilon
G^+(\Lambda _\pm )\psi $  belong to $\frak S_{2(n-1)/3}\cup\Cal C_1$. Similarly, $
G^-(\Lambda _\pm)\zeta _\varepsilon $ are of order $-\frac 32$, and 
$G^-(\Lambda _\pm )\zeta _\varepsilon \psi $,  $\psi G^-(\Lambda _\pm )\zeta _\varepsilon $  
belong to $\frak S_{2(n-1)/3}\,\cup\,\Cal C_1$. 
\endproclaim 

\demo{Proof} It suffices to give the details for $\varepsilon
=1$. Consider $G^+(\Lambda _+)$.
First we note that
$$
\zeta G^+(\Lambda _+)=
\zeta r^+\Lambda _+ e^-J=r^+\Lambda _+ \zeta e^-J+ r^+[\zeta ,\Lambda _+
]e^-J=
r^+[\zeta ,\Lambda _+ ]e^-J,
$$
since $\zeta e^-=0$; here $[\zeta ,\Lambda _+]$ is the commutator $\zeta \Lambda _+-\Lambda
_+\zeta $. As for ordinary $\psi $do's, the commutator is of lower
order; since $\Lambda _+$ is nonstandard, we work out proof details:

For
$t,s\in{\Bbb R}$, $\zeta $ has the Taylor-expansion
$$
\aligned
\zeta (t)&=\sum_{0\le j<J}\tfrac1{j!}\zeta
^{(j)}(s)(t-s)^j+(t-s)^J\varrho _J(s,t), \text{ where }\\
\varrho _J(s,t)&=\tfrac 1{(J-1)!}\int_0^1(1-h)^{J-1}\partial^J\zeta
(s+h(t-s))\, dh.
\endaligned
$$ Then, using that $(t-s)^je^{i(t-s)\tau }=D_\tau ^je^{i(t-s)\tau }$
and integrating by parts (as allowed in oscillatory integrals), we
find,  denoting $(2\pi )^{1-n}d\xi '=\d \xi '$:
$$\aligned
[\zeta ,\Lambda _+]u&=\int e^{i(x'-y')\cdot\xi '}(\zeta (t)-\zeta
(s))\lambda _+(\xi ')u(y')\, \d \xi 'dy'\\
&=\int e^{i(x'-y')\cdot\xi '}(\sum_{1\le j<J}\tfrac1{j!}\zeta
^{(j)}(s)(t-s)^j+(t-s)^J\varrho _J (s,t))\lambda _+(\xi ')u(y')\,
\d \xi 'dy'\\
&=\int e^{i(x'-y')\cdot\xi '}(\sum_{1\le j<J}\tfrac1{j!}\zeta
^{(j)}(s)\overline D_\tau ^j+\varrho _J (s,t)\overline D_\tau ^J)\lambda
_+(\xi ')u(y')\, \d \xi 'dy'\\
&=\sum_{1\le j<J}\tfrac1{j!}\Op(\overline D_\tau ^j\lambda _+(\xi '))
\zeta
^{(j)} u+\Op(\varrho _J (s,t)\overline D_\tau ^J\lambda
_+(\xi '))u.
\endaligned\tag5.13$$
Here
$$
\overline D_\tau ^j\lambda _+(\xi ')=\overline D_\tau ^j(|(\xi '',\alpha 
)|+i\tau )^{-\frac12}=c_j (|(\xi '',\alpha 
)|+i\tau )^{-\frac12-j};
$$
they are of order $-\frac12-j$. Take $J$ so large
that the last symbol is integrable in $\xi '$, e.g.\ $J=n$. Then the
terms in the sum over $j$ map $H^r({\Bbb R}^{n-1})$ into
$H^{r+\frac12+j}({\Bbb R}^{n-1})$ ($r\in{\Bbb R}$)
by elementary considerations, and the last term has a continuous
kernel, supported for $s,t\in [\frac13,\frac43]$. Similar
considerations hold for $\zeta G^+(\Lambda _-)$.
When we cut down with multiplication by $\psi $, and functions
$1_{{\Bbb R}^{n-1}_\pm}$, we can use the spectral
estimates (5.4) and the trace-class property of
operators with continuous kernel, to see that $\psi \zeta G^+(\Lambda
_\pm)$ are in $\frak S_{2(n-1)/3}\,\cup\,\Cal C_1$. .

The statements for $G^-(\Lambda _\pm)\zeta $ are shown similarly, and for
the operators with $\psi $ to the right one can
use that $G^+(\Lambda _\pm)$ and $G^-(\Lambda _\mp)$ are adjoints. \qed 
\enddemo

One could argue in a more refined way (e.g.\ with sequences of nested
cutoff functions), to show that since $\zeta $ is
supported away from 0, {\it all} the terms in
(5.13) give spectrally negligible contributions (in $\bigcap_p\frak S_p$) when we take $G^+$ of
them 
(as for ordinary singular Green operators), but that extra
information will not be needed in the following.

Next, we shall show spectral estimates of the contributions to
$G^\pm(\Lambda _\pm)$ supported near $t=0$. Here we shall profit from
the fact that Birman and Solomyak in \cite{BS77} showed far-reaching
spectral results for nonstandard $\psi $do's, taking $L_p$-norms (not just
$L_\infty $-norms) of cutoff functions into account. Anisotropic
symbols are allowed there, but we just need the case of isotropic
symbols with low smoothness.

\proclaim{Theorem 5.4}{\rm \cite{BS77}} Let $A=\Op(b(x)a(x,\xi
)c(y))$ on ${\Bbb R}^m$, with $a(x,\xi )$ homogeneous in $\xi $ of
degree 
$-\mu  \in
\,]-m,0[\,$. Denote $m/\mu  =\nu $. Then $A\in \frak S_\nu $ with  
$$
\sup _{j\in{\Bbb N}}s_j(A)j^{1/\nu }\le C \|b\|_{L_{q_1}}\|c\|_{L_{q_2}}[a_{|\xi |=1}]_\beta ,
$$
 if 
$$
q_1,q_2\in \,]2,\infty ],\quad \tfrac1{q_1}+\tfrac1{q_2}=\tfrac1{\nu
},\quad \beta =q_1. 
$$
(Here $[\Phi (x,\xi ) ]_\beta $ denotes the norm of a certain linear operator
on $\frak S_\beta $ defined from $\Phi $.) A sufficient condition for
the boundedness of $[a_{|\xi |=1}]_\beta $ is that
$$
a(x,\xi )|_{|\xi |=1}\in L_\infty (S_\xi ^{m-1}, W^l_p({\Bbb
R}^m_x)),\text{ with }\tfrac12-\tfrac1{q_1}<\tfrac1p\le \tfrac12,\;
 p\,l>m. \tag5.14
$$
\endproclaim

The paper \cite{BS77} also covers cases where $\nu \le 1$, and gives
spectral asymptotics formulas under additional mild regularity
hypotheses (in (5.14), $L_\infty $ is then replaced by $C^0$). 

In order to apply the result we must estimate the effect of replacing
the $\Lambda _{\pm}(\alpha )$ by the operators $\Lambda _{\pm}(0)$ with
strictly homogeneous symbols $\lambda _{\pm}(\xi ',0)$.

\proclaim{Lemma 5.5} The symbol $\lambda _{\pm}(\xi ',\alpha )-\lambda
_{\pm}(\xi ',0)$ of the difference $\Lambda _{\pm}(\alpha )-\Lambda
_{\pm}(0)$ satisfies$$
\lambda _{\pm}(\xi ',\alpha )-\lambda
_{\pm}(\xi ',0)=
O(\alpha ^2|(\xi ',\alpha )|^{-\frac32}|\xi '|^{-1}|(\xi
'',\alpha )|^{-1}).
$$ Hence it defines an operator mapping $H^r({\Bbb R}^{n-1})$ into
$H^{r+\frac52}_{\operatorname{loc}}({\Bbb R}^{n-1})$ for $r\in{\Bbb
R}$. 

It follows that $\psi (\Lambda _+(\alpha )-\Lambda _+(0))$, $\psi
(G^+(\Lambda _+(\alpha ))-G^+(\Lambda _+(0)))$, $ (\Lambda _-(\alpha
)-\Lambda _-(0))\psi _1$ and $ (G^-(\Lambda _-(\alpha ))-G^-(\Lambda _-(0)))\psi _1$
are in $\frak S_{2(n-1)/5}$. 
\endproclaim

\demo{Proof} We give the details for $\Lambda _+(\alpha )-\Lambda _+(0)$.
Here
$$
\aligned
\lambda _+(\xi ',\alpha ) -\lambda _+(\xi ',0 ) &=\tfrac{\lambda _+(\xi
',\alpha )^2- \lambda _+(\xi ',0 )^2}{\lambda _+(\xi ',\alpha  )+\lambda
_+(\xi ',0 )}
=\tfrac1{\lambda _+(\xi ',\alpha  )+\lambda _+(\xi ',0
)}\big(\tfrac1{|(\xi '',\alpha )|+i\tau }-\tfrac1{|\xi ''|+i\tau }\big)\\
&=\tfrac{|\xi ''|-|(\xi '',\alpha )|}{(\lambda _+(\xi ',\alpha  )+\lambda
_+(\xi ',0 ))(|(\xi '',\alpha )|+i\tau )(|\xi ''|+i\tau )}\\
&=\tfrac{-\alpha ^2}{(\lambda _+(\xi ',\alpha  )+\lambda _+(\xi ',0 ))(|(\xi
'',\alpha )|+i\tau )(|\xi ''|+i\tau )(|\xi ''|+|(\xi '',\alpha )|)}\\
&=O(\alpha ^2|(\xi ',\alpha )|^{-\frac32}|\xi '|^{-1}|(\xi '',\alpha )|^{-1}).
\endaligned
$$
The operator with symbol $\zeta (|\xi '|)(\lambda _+(\xi ',\alpha ) -\lambda _+(\xi ',0 ))$
maps $H^r({\Bbb R}^{n-1})$ into
$H^{r+\frac52}({\Bbb R}^{n-1})$ for $r\in{\Bbb R}$, and the remainder
supported near $|\xi '|=0$ gives an operator mapping into $C^\infty
({\Bbb R}^{n-1})$. When cutoffs by compactly supported functions are
 applied, this gives operators in $\frak S_{2(n-1)/5}$.

The result for $\Lambda _-$ follows by similar calculations or by duality. 
\qed
\enddemo

In the following, $\varphi (t)$ denotes a function in $C^\infty ({\Bbb
R})$ that  takes
values in $[0,1]$ and equals $1$ for $|t|\le \frac13$, equals $0$ for
$|t|\ge \frac23$; we denote $\varphi  (t/\varepsilon
)=\varphi  _\varepsilon (t)$. We can assume that $1_{\rp}(1-\zeta
)=1_{\rp}\varphi $.

\proclaim{Lemma 5.6} There are the following spectral estimates:
$$
\aligned
\sup _js_j(\varphi _\varepsilon (t )\psi \Lambda _+(0))j^{1/(2n-2)}&\le
C_\varepsilon ,\\
\sup _js_j(\varphi _\varepsilon (t )\psi G^+(\Lambda _+(0)))j^{1/(2n-2)}&\le
C_\varepsilon ,\\
\sup _js_j( \Lambda _-(0)\psi _1\varphi _\varepsilon (t ))j^{1/(2n-2)}&\le
C_\varepsilon ,\\
\sup _js_j(G^-( \Lambda _-(0))\psi _1\varphi _\varepsilon (t ))j^{1/(2n-2)}&\le
C_\varepsilon .
\endaligned\tag5.15
$$
where $C_\varepsilon \to 0 $ for $\varepsilon \to 0$.

\endproclaim

\demo{Proof} 
For the first line in (5.15), we apply Theorem 5.4 with
$$
b(x')=\varphi _\varepsilon (t)\psi _2(x'),\quad a(x',\xi ')=\psi
(x')\lambda _+(\xi ',0), \quad c(x')=1,
$$
where $\psi _2\in C_0^\infty ({\Bbb R}^{n-1})$, equal to 1 on $\supp \psi $.
Here 
$m=n-1$, $\mu  =\frac12$ so that $\nu =m/\mu =2(n-1)$, and we take
$q_1=\beta =\nu =2(n-1)$ and $q_2=\infty $. Moreover, since
$\frac12-\frac1{q_1}=\frac12-\frac1{2(n-1)}=\frac{n-2}{2(n-1)}$, 
$p$ is taken in $\,]2,\frac {2(n-1)}{n-2}]$ and $l$ is taken
$>(n-1)/p$. (5.14) is satisfied since $\psi \in C_0^\infty ({\Bbb
R}^{n-1})\subset W^l_p({\Bbb R}^{n-1})$. Then
$$\aligned
\sup _js_j(\varphi _\varepsilon (t )\psi \Lambda
_+)j^{1/(2n-2)}&\le C\|\varphi _\varepsilon \psi _2\|_{L_{\nu }}
\le C'\operatorname{vol}(\supp (\varphi _\varepsilon \psi
_2))^{1/(2n-2)}\\
&\le
C''\varepsilon ^{1/(2n-2)}\to 0
\endaligned$$
for $\varepsilon \to 0$.

For the second line in (5.15) we replace $b$ by $1_{{\Bbb
R}^{n-1}_+}\varphi _\varepsilon \psi _2$ and $c$ by $1_{{\Bbb
R}^{n-1}_-}$, and use that $J$ is an isometric isomorphism.

The proof of the third and fourth line goes in a similar way, 
interchanging choices for $b$
and $c$.
\qed
\enddemo

We can finally conclude:

\proclaim{Theorem 5.7} The operator $L_0^{-1}$ acts like 
$$
L_0^{-1}= -P_{0,\nu ,\gamma ,+}-G^+(\Lambda _+)G^{-}(\Lambda _-),\tag5.16
$$
where 
$$
P_{0,\nu ,\gamma ,+}\colon  H^{s }(\rp^{n-1})\to
H^{s+1 }(\rp^{n-1})\text{ for }-\tfrac12<s<\tfrac12,\tag5.17
$$
and the operators $\psi G^{\pm}(\Lambda
_{\pm})$ and $G^{\pm}(\Lambda _{\pm})\psi $ are in  $\frak S_{2(n-1),0}$,
when
$\psi \in C_0^\infty ({\Bbb R}^{n-1})$.\endproclaim

\demo{Proof} 
The decomposition (5.16) was shown above. The continuity in (5.17)
follows since $P_{0,\nu ,\gamma }$ is a constant-coefficient $\psi $do of order $-1$.
For the next statement, we give details for $\psi G^+(\Lambda _+)$; the other
cases are similar. For any $\varepsilon >0$ we can write
$$
\psi G^+(\Lambda _+(\alpha ))
=\varphi _\varepsilon (t)\psi G^+(\Lambda
_+(0))+\varphi _\varepsilon (t)\psi [G^+(\Lambda _+(\alpha ))-G^+(\Lambda
_+(0))]+\zeta _\varepsilon (t)
\psi G^+(\Lambda _+(\alpha )).\tag5.18
$$
Here the first term satisfies (5.15), the second term is in $\frak
S_{2(n-1)/5}$ by Lemma 5.5, and the third term is in $\frak S_{\max
\{2(n-1)/3, 1+\delta \}}$ (for any $\delta >0$) by Lemma 5.3. Thus the sum
of the second and third term satisfies $s_j j^{1/(2n-2)}\to 0$ for
$j\to\infty $.
We can then apply Lemma 5.1 $2^\circ$, with $1/p= 1/(2n-2)$,
$M=1/\varepsilon $, $B_M$ being the sum of the second and third terms
and $B'_M$ being the first term, $c_M=C_\varepsilon $ and $C_M=C_0=0$.\qed 
\enddemo

\example{Remark 5.8} In the case $n=2$, when
$Q$ is a $\psi $do on ${\Bbb R}^{n-1}$ of order $-r<0$,  the operators $G^\pm(Q)$ are not covered by
Theorem 5.2. But certainly the calculations leading to Theorem 5.7
work in this case, so we have at least that
$\psi G^\pm(Q)\in\frak S_{(n-1)/r,0}$.
Similarly, if $n=2$ and $Q_1$ and $Q_2$ are $\psi $do's on $\Sigma $
of negative orders $-r_1,-r_2$, then 
$L(Q_1,Q_2)\in \frak S_{(n-1)/(r_1+r_2)}$.
 \endexample

\subhead 5.3 Variable coefficients, analysis of $L^{-1}$ \endsubhead

Now consider $A=-\Delta +a_0(x)$ on the smooth bounded open subset
$\Omega $ of ${\Bbb R}^n$, provided with the mixed boundary condition
$\nu u=b\gamma _0u$ on $\Sigma _+$, $\gamma _0u=0$ on $\Sigma _-$. The
operator $L$ acts like 
$$
L\varphi =r^+(b-P_{\gamma ,\nu })e^+\varphi =-P_{\gamma ,\chi ,+}\varphi 
$$
for $\varphi \in D(L)$, cf.\ (4.3), (4.6).

In the analysis of $L^{-1}$ on $\Sigma _+$ we want to use the insight
gained in Section 5.2 for the ``flat'' constant-coefficient case, but since the
ingredients are not standard $\psi $do's, we do not have the usual
localization tools for $\psi $do's available and must reason very
carefully (for example in Eskin's book, formulas for coordinate
changes are only worked out for a subclass of symbols with better
estimates than the present $\lambda _{\pm}(\xi ')$).
The strategy will be to reduce to a situation where the results from
the   ``flat'' case can be used directly. 

Our aim is to show:

\proclaim{Theorem 5.9} The operator $L^{-1}$ acts like $-P_{\nu ,\gamma ,+}+R$,
where $R\in \frak S_{n-1,0}$. In particular, $L^{-1}\in \frak S_{n-1}$.
\endproclaim

This will be shown in several steps. We first show a preliminary
spectral estimate for $L^{-1}$; it will be improved later.

\proclaim{Lemma 5.10} The operator $L^{-1}\colon X^*\to X$ extends to an operator $M$ that maps continuously
$$
M\colon H^s(\Sigma _+^\circ)\to H^{s+\frac12-\varepsilon }_0(\Sigma
_+)\text{ for }-1<s\le \tfrac12.\tag5.19
$$
In particular, the closure of $L^{-1}$ in $L_2(\Sigma _+)$ is a continuous operator 
from $L_2(\Sigma _+)$ to $H^{\frac12-\varepsilon }_0(\Sigma _+)$;
it  belongs to $\frak S_{(n-1)/(\frac12-\varepsilon )}$ for $\varepsilon >0$.
\endproclaim  

\demo{Proof} It follows from Proposition 4.1 that $L^{-1}$ is
continuous from $X^*=H^\frac12(\Sigma _+^\circ)$ to $H^{1-\varepsilon
}_0(\Sigma _+)$. Then it has an adjoint $M$ (with respect to dualities
consistent with the $L_2(\Sigma _+)$-scalar product) that is
continuous from $H^{-1+\varepsilon }(\Sigma _+^\circ)$ to
$H^{-\frac12}_0(\Sigma _+)$. But since $L^{-1}$ is known to be
selfadjoint (from $X^*$ to $X$, consistently with the $L_2$-scalar product), $M$ must be an extension of
$L^{-1}$. Now (5.19) follows  by interpolation. For $s=0$ we find the
last statement, where the spectral information follows from (5.4);
note that $H^{\frac12-\varepsilon }_0(\Sigma
_+)=H^{\frac12-\varepsilon }(\Sigma _+^\circ)$.\qed 
\enddemo

When $\{\varrho
_1,\dots,\varrho _N\}$ is any partition of unity for $\Sigma $, then
$L^{-1}=\sum_{k=1}^N\varrho _kL^{-1}$, and it suffices to analyze the terms
$\varrho _kL^{-1}$ individually. Here we  can also introduce a cutoff
function $\psi _k$ to the right, considering terms $\varrho _kL^{-1}\psi _k$
where $\psi _k$ is 1 on the support of $\varrho _k$; the effect of
such a modification
will be studied later. 

Our next observation is that it is allowed to perform smooth diffeomorphisms
of $\Omega $, in particular of $\Sigma $. Assume that $\kappa $ is a
diffeomorphism of an open neighborhood $U_0$ of $\comega$ onto another
open set $V_0\subset {\Bbb R}^n$, where $\kappa (\Omega )=\Omega '$,
then functions $f(x)$ on $\Omega $ are carried over to functions $\underline
f(y)=f(\kappa ^{-1}(y))$ on $\Omega '$, and operators $P$ over $\comega$ are carried
over to operators $\underline P$ over $\overline{\Omega'}$:
$$
(\underline P\underline f)(y)=(Pf)(\kappa ^{-1}(y)).\tag5.20
$$
   The $\psi $do $P_{\gamma ,\chi }$ on $\Sigma $ carries over to a
   $\psi $do $\underline {P_{\gamma ,\chi }}$ on $\underline\Sigma $ according to well-known rules; it
   is again elliptic of order 1 and has the same principal symbol. The
   operator $L$ carries over to $\underline L$, equal to the truncated
   version of  $\underline {P_{\gamma ,\chi }}$, where we apply $e^+$
   and $r^+$ with respect to the partition $\underline\Sigma
   =\underline\Sigma _+\cup \underline\Sigma _-$. There is again an
   inverse $\underline L^{-1}$, with mapping properties as explained
   for $L^{-1}$, relative to the transformed sets.

To find the structure of $L^{-1}$ in a neighborhood of a point
$x_0\in\Sigma $, let us consider $\psi L^{-1}\psi _1$, where $\psi $
and $\psi _1$ are $C^\infty $-functions supported in the
neighborhood, with $\psi _1=1 $ on $\operatorname{supp}\psi $.

It can be assumed, after a translation and rotation if necessary, that
$x_0\in \crp^{n-1}$ and the interior normal at $x_0=\{x_{0,1},\dots, x_{0,n-1},0\}$ is
$(0,\dots,0,1)$, such that $x_{0,n-1}>0$ if $x_0\in \Sigma _+^\circ$
and $x_{0}=0$ if $x_0\in \partial\Sigma _+$; in the latter case we can
assume that the interior normal to
$\partial\Sigma _+\subset \Sigma $ at $x_0$ is $\{0,\dots,0,1,0\}$.
We choose a diffeomorphism that changes $\comega$ only near $x_0$.
If $x_0\in \Sigma _+^\circ$, we can assume that $\psi $ and $\psi _1$
are supported away from
$\partial\Sigma _+$; then we let the diffeomorphism be such that it  
transforms
a neighborhood
$U\subset{\Bbb R}^n$ of $x_0$ over to  $V\subset{\Bbb R}^n$, carrying
$U\cap\Omega $ and $U\cap \Sigma _+$ over to
$V\cap \rnp$ and $V\cap \crp^{n-1}$, with $\psi $
and $\psi _1$ supported in $U\cap \Sigma _+^\circ$.
If $x_0\in \partial\Sigma _+$, we choose the diffeomorphism such that
$U\cap\Omega $, $U\cap \Sigma $ and $U\cap \Sigma _+$ are mapped to
$V\cap {\Bbb R}^np$, $V\cap {\Bbb R}^{n-1}$ and $V\cap
\crp^{n-1}$, $\psi $ and $\psi _1$ supported in $V\cap {\Bbb R}^{n-1}$. 
(The identifications of ${\Bbb R}^{n-1}$ and $\crp^{n-1}$ with  ${\Bbb
R}^{n-1}\times\{0\}$ and $\crp^{n-1}\times\{0\}$ as subsets of ${\Bbb
R}^n$ are understood here.)  

This gives a transformed operator $\underline \psi \underline
L^{-1}\underline\psi _1$ acting on functions supported in
$V'=V\cap{\Bbb R}^{n-1}$. For simplicity of notation, we drop the
underlines in the following.

We shall compare $\psi L^{-1}\psi _1$ with $\psi L_0^{-1}\psi _1$
where $L_0^{-1}$ is the constant-coefficient operator studied in
Section 5.2. Let us give the details for the most delicate case
$x_0\in \partial\Sigma _+$, where the effects of truncation have to be
taken into account.

\proclaim{Proposition 5.11} In the setting described in the preceding lines, we have that 
$$
\psi L^{-1}\psi _1=-\psi P_{\nu ,\gamma ,+}\psi _1+R_1,\tag5.21
$$
as operators in $L_2(V'\cap \rp^{n-1})$, where $R_1\in \frak S_{n-1,0}$.
\endproclaim

\demo{Proof} It follows from Theorem 5.7 that
$$
\psi L_0^{-1}\psi _1=-\psi P_{0,\nu ,\gamma ,+}\psi _1+R_2,\tag5.22
$$
where $R_{2}=\psi G^+(\Lambda _+)G^-(\Lambda _-)\psi _1$ is in $\frak
S_{n-1,0}$; cf.\ (5.2). We shall now compare $\psi L^{-1}\psi _1$ and
$\psi L_0^{-1}\psi _1$. There is the difficulty that the operators
$L_1^{-1}$ and $L_0^{-1}$ do not act over the same manifold, but this will
be dealt with by introduction of more cutoff functions.
Let $\psi _2\in C_0^\infty (V')$, satisfying $\psi _2=1$ on $\supp
\psi _1$. We calculate:
$$
\aligned
\psi L^{-1}\psi _1-\psi L_0^{-1}\psi _1&=\psi L^{-1}\psi _2\psi
_1-\psi \psi _2L_0^{-1}\psi _1\\
&=\psi L^{-1}\psi _2L_0L_0^{-1}\psi _1-\psi L^{-1}L\psi _2L_0^{-1}\psi _1.
\endaligned \tag5.23$$
We want to insert the factor $\psi _2$ in the middle of $L_0L_0^{-1}$
as well as $L^{-1}L$; this is justified as follows:
Write e.g.
$$
\psi L^{-1}\psi _2L_0L_0^{-1}\psi _1=\psi L^{-1}\psi _2L_0\psi _2L_0^{-1}\psi _1+\psi L^{-1}\psi _2L_0(1-\psi _2)L_0^{-1}\psi _1.
$$
For the last term, we note that (since $(1-\psi _2)\psi _1=0$)
$$
(1-\psi _2)L_0^{-1}\psi _1=[(1-\psi _2),L_0^{-1}]\psi _1=[L_0^{-1},\psi _2]\psi _1=L_0^{-1}[\psi _2,L_0]L_0^{-1}\psi _1,
$$
where $[L_0,\psi _2]=[-P_{0,\gamma ,\nu ,+},\psi _2]$ is $L_2$-bounded (since $P_{0,\gamma ,\nu }$ is a first-order $\psi $do). Then
$$
\psi L^{-1}\psi _2L_0(1-\psi _2)L_0^{-1}\psi _1=\psi L^{-1}\psi _2L_0L_0^{-1}[\psi _2,L_0]L_0^{-1}\psi _1=\psi L^{-1}\psi _2[\psi _2,L_0]L_0^{-1}\psi _1,
$$
which is the composition of $\psi L^{-1}\psi _2\in \frak
S_{(n-1)/(\frac12-\varepsilon )}$ (cf.\ Lemma 5.10), the bounded
operator $[\psi _2,L_0]$, and
$L_0^{-1}\psi _1\in \frak S_{n-1}$ (cf.\ Theorem
5.7; its adjoint is $\overline \psi _1L_0^{-1}$). Then the whole term is in
$\frak S_{(n-1)/(\frac32-\varepsilon )}$, and 
$$
\psi L^{-1}\psi _2L_0L_0^{-1}\psi _1=\psi L^{-1}\psi _2L_0\psi
_2L_0^{-1}\psi _1+R_3,\text{ where }R_3\in \frak S_{(n-1)/(\frac32-\varepsilon )}.
$$

Similarly, we can insert a factor $\psi _2$ between $L^{-1}$ and $L$
in the last term of (5.23), making an error that is in $\frak
S_{(n-1)/(\frac32-\varepsilon )}$. 

It remains to consider
$$
\psi L^{-1}\psi _2L_0\psi _2L_0^{-1}\psi _1-\psi L^{-1}\psi _2L\psi
_2L_0^{-1}\psi _1=
(\psi L^{-1}\psi _3)(\psi _2L_0\psi _2-\psi _2L\psi _2)(\psi _3L_0^{-1}\psi _1),
$$
 where we have replaced $\psi _2$ by $\psi _2\psi _3$, with $\psi _3=1$
 on $\supp \psi _2$, in a few places. Here the first factor is in $\frak
 S_{(n-1)/(\frac12-\varepsilon )}$ by Lemma 5.10, the last factor is in $\frak
 S_{n-1}$ by Theorem 5.7, and the middle factor is a truncated $\psi $do of order
 zero, hence bounded in $L_2$, since $P_{\gamma ,\chi  }$ and $P_{0,\gamma ,\nu }$ have the
 same principal symbol on $V'$ (recall (4.6)). Then the whole expression is in $\frak
 S_{(n-1)/(\frac32-\varepsilon )}$. Thus we have obtained that
$\psi L^{-1}\psi _1-\psi L_0^{-1}\psi _1\in \frak
 S_{(n-1)/(\frac32-\varepsilon )}$, which is contained in 
$ \frak S_{n-1,0}$. 
Together with  (5.22) this shows
$$
\psi L^{-1}\psi _1-\psi (-P_{0,\nu ,\gamma ,+})\psi _1\in \frak S_{n-1,0}. 
$$

Finally, since $P_{0,\nu ,\gamma }$ and $P_{\nu ,\gamma }$ have the same
principal symbol (of order $-1$) on $V'$, $\psi P_{0,\nu ,\gamma
,+}\psi _1-\psi P_{\nu ,\gamma ,+}\psi _1$ is a truncated $\psi $do of
order $-2$; hence it is in $\frak S_{(n-1)/2}\subset \frak S_{n-1,0}$,
and (5.21) follows.\qed 
\enddemo

\demo{Proof of Theorem 5.9} We now consider $L^{-1}$ on $\Sigma
_+$, written as  $L^{-1}=\sum_{k=1}^N\varrho _kL^{-1}$ for some partition of
unity $\sum_{k=1}^N\varrho _k=1$. To analyze an individual term
$\varrho _kL^{-1}$, we choose a cutoff function $\psi _k$ that is 1
on $\supp \varrho _k$, and write$$\aligned
\varrho _kL^{-1}&=\varrho _kL^{-1}\psi _k+\varrho _kL^{-1}(1-\psi
_k)=\varrho _kL^{-1}\psi _k+\varrho _k[L^{-1},1-\psi _k]\\
&=\varrho _kL^{-1}\psi _k+\varrho _kL^{-1}[L,\psi _k]L^{-1}.
\endaligned\tag5.24$$ 
Since $[L,\psi _k]$ is a truncated zero-order $\psi $do, it is bounded
in $L_2$. By Lemma 5.10, $L^{-1}\in \frak
S_{(n-1)/(\frac12-\varepsilon )}$, so the last term satisfies
$$
\varrho _kL^{-1}[L,\psi _k]L^{-1}\in \frak S_{n-1+\delta }, \text{ any
}\delta >0.\tag5.25
$$

We can assume that the supports of $\varrho _k$ and $\psi _k$ are so
small that a diffeomorphism as described before Proposition 5.11 can
be applied in a neighborhood of the supports; then Proposition 5.11
gives that
$$
\varrho _k L^{-1}\psi _k=-\varrho _kP_{\nu ,\gamma ,+}\psi
_k+R_{1,k},\text{ with }R_{1,k}\in \frak S_{n-1,0}.\tag5.26
$$
A first observation resulting from this is that $\varrho _k L^{-1}\psi
_k$ is in $\frak S_{n-1}$, since the $\psi $do of order $-1$ is
there. Then in view of (5.24)--(5.25), $\varrho _kL^{-1}\in \frak
S_{n-1+\delta }$, any $\delta >0$. Summation in $k$ gives that  $L^{-1}\in \frak
S_{n-1+\delta }$. Next, we go back to (5.24), where the new
information allows us to conclude that
$$
\varrho _kL^{-1}[L,\psi _k]L^{-1}\in \frak S_{(n-1+\delta )/2}\subset
\frak S_{n-1,0}.
$$
In view of (5.26), we finally get that 
$$
\varrho _kL^{-1}=-\varrho _kP_{\nu ,\gamma ,+}\psi
_k+R_{2,k},\text{ with }R_{2,k}\in \frak S_{n-1,0}.
$$
Summation in $k$ gives that
$$
L^{-1}=-\sum_{k=1}^N\varrho _kP_{\nu ,\gamma ,+}\psi
_k+\sum _{k=1}^NR_{2,k}=-P_{\nu ,\gamma ,+}+R_{3},
$$
with $R_{3}\in \frak S_{n-1,0}$. \qed

\enddemo

\subhead 5.4 Reduction of the Poisson operators \endsubhead

We now consider the operator (5.1). 
To find the spectral behavior, we note that by the general
rule for eigenvalues $\mu _j(ST)=\mu _j(TS)$, we can write
$$
\mu _{j}(\inj_{V  }\gamma _V^{-1}L
^{-1}(\gamma _V^{-1})^*\pr_{V  })
=\mu _{j}(L
^{-1}(\gamma _V^{-1})^*\pr_{V  }\inj_{V  }\gamma _V^{-1})
=\mu _{j}(L
^{-1}(\gamma _V^{-1})^*\gamma _V^{-1}),\tag5.27
$$
in view of (2.15). 

\proclaim{Lemma 5.12} The operator $(\gamma _V^{-1})^*\gamma _V^{-1}$
satisfies
$$
(\gamma _V^{-1})^*\gamma _V^{-1}=P_{1,+},\tag5.28
$$
where $P_1=K_\gamma ^*K_\gamma $ is a selfadjoint nonnegative elliptic $\psi $do of
order $-1$ on $\Sigma $ with principal symbol $(2|\xi '|)^{-1}$.
\endproclaim

\demo{Proof}
We have for $\varphi,\psi  \in X$:
$$
\multline
((\gamma _V^{-1})^*\gamma _V^{-1}\varphi ,\psi )_{X^*,X}=
(\gamma _V^{-1}\varphi ,\gamma _V^{-1}\psi )_{V}=
(K_\gamma \varphi ,K_\gamma \psi )_H\\
=(K_\gamma ^*K_\gamma \varphi ,\psi )_{\frac12,-\frac12}=(P_1\varphi
,\psi )_{\frac12,-\frac12},\endmultline \tag5.29
$$
where 
$P_1=K_\gamma ^*K_\gamma $ is a $\psi $do of
order $-1$ on $\Sigma $, by the rules of calculus for
pseudodifferential boundary operators; it is clearly selfadjoint nonnegative. The principal symbol is found from the calculation using  the
principal symbol-kernel $\tilde k^0=e^{-x_n|\xi '|}$ of $K_\gamma $:
$$
\int_0^\infty e^{-x_n|\xi '|}e^{-x_n|\xi '|}\,dx_n=(2|\xi '|)^{-1},
$$
also equal to $\|\tilde k^0\|^2_{L_2(\rp)}$.
Since $\varphi $ and $\psi $ are supported in
$\Sigma _+$, (5.29) may be rewritten further as
$$
(P_1\varphi ,\psi )_{\frac12,-\frac12}=(P_1e^+\varphi ,e^+\psi
)_{\frac12,-\frac12}
=(r^+P_1e^+\varphi ,\psi )_{H^\frac12(\Sigma ^\circ_+), H^{-\frac12}_0(\Sigma _+)}
=(P_{1,+}\varphi , \psi )_{X^*,X}.
$$
Then (5.28) follows since $\varphi $ and $\psi $ are arbitrary.\qed
\enddemo

Next, we define 
$$
P_2=P_1^{\frac12},\tag5.30
$$
a nonnegative selfadjoint $\psi $do on $\Sigma $ of order $-\frac12$,
by Seeley \cite{S67}. Moreover, set
$$
G^{(1)}=P_{1,+}-(P_{2,+})^2,\quad G^{(\frac12)}=(P_{1,+})^{\frac12}-P_{2,+}.\tag5.31
$$

\proclaim{Lemma 5.13} When $n\ge 3$, $G^{(1)}\in \frak S_{n-2}$ and $G^{(\frac12)}\in
\frak S_{2(n-2)}$. When $n=2$, $G^{(1)}\in \frak S_{n-1,0}$ and $G^{(\frac12)}\in
\frak S_{2(n-1),0}$.
\endproclaim
\demo{Proof} We first note (cf.\ (5.12)) that
$$
G^{(1)}=P_{1,+}-P_{2,+}P_{2,+}=r^+P_2P_2e^+-r^+P_2e^+r^+P_2e^+=L(P_2,P_2).\tag5.32
$$
Since $P_2$ is a $\psi $do of order $-\frac12$, we have by Theorem
5.2 that $L(P_2,P_2)$
is in $\frak S_{n-2}$ when $n\ge 3$. For $n=2$, we see that
$s_j(L(P_2,P_2))j^{1/(n-1)}\to 0$ for $j\to\infty $ by use of Remark 5.8. 

To obtain the result for $G^{(\frac12)}$, we shall as in \cite{G83}
appeal to a result of Birman, Koplienko and Solomyak \cite{BKS75}. It states that
when $M_1$ and $M_2$ are compact selfadjoint nonnegative operators on a Hilbert
space $H$ such that $G^{(1)}=M_1-M_2$ is in $\frak S_\gamma $ for some
$\gamma >0$, then  $G^{(\sigma )}=M_1^\sigma -M_2^\sigma $ is in
$\frak S_{\gamma/\sigma } $ for all $0<\sigma <1$. Applying this with
$M_1=P_{1,+}$, $M_2=(P_{2,+})^2$ and $\sigma =\frac12$, we get the
desired result when $n\ge 3$. The paper \cite{BKS75} also shows that
$\lim\sup s_j(G^{(\frac12)})j^{1/(2n-2)}$  is dominated by $\lim\sup
s_j(G^{(1)})j^{1/(n-1)}$, which assures the statement
for $n=2$.\qed  
\enddemo

Now we continue the analysis in (5.27) as follows:

\proclaim{Proposition 5.14}
$$
\mu _{j}(\inj_{V  }\gamma _V^{-1}L
^{-1}(\gamma _V^{-1})^*\pr_{V  })=\mu _{j}(P_{2,+}
L
^{-1}P_{2,+} +G'),\tag5.33
$$
where $G'$ is the selfadjoint operator
$$
G'=
G^{(\frac12)}L^{-1}P_{2,+}+P_{2,+}L^{-1}G^{(\frac12)}+G^{(\frac12)}L^{-1}G^{(\frac12)}
;\tag5.34
$$
it is in $\frak
S_{(n-1)/2-r}$ for a positive $r$ when $n\ge 3$, and in $\frak
S_{(n-1)/2,0}$ when $n=2$.
\endproclaim 

\demo{Proof}
Using Lemma 5.13 and (5.27) and the definitions (5.30)--(5.31) we have:
$$\multline
\mu _{j}(\inj_{V  }\gamma _V^{-1}L
^{-1}(\gamma _V^{-1})^*\pr_{V  })=\mu _{j}(L
^{-1}P_{1,+})= \mu _{j}(L
^{-1}(P_{1,+})^\frac12(P_{1,+})^\frac12)\\= \mu _{j}((P_{1,+})^\frac12
L^{-1}(P_{1,+})^\frac12)=
\mu _{j}(P_{2,+}
L
^{-1}P_{2,+} +G'),\endmultline
$$
where $G'$ is as in (5.34).
When $n\ge 3$, we use that $L^{-1}\in \frak S_{n-1 }$, $P_{2,+}\in \frak S_{2(n-1)}$,
and $G^{(\frac12)}\in \frak S_{2(n-2)}$ ( by Lemma 5.13) and the rule
(5.2) to see that$$
G'\in\frak S_p\text{ with }p=(\tfrac1{n-1}+\tfrac
1{2(n-1)}+\tfrac 1{2(n-2)})^{-1}<\tfrac{n-1}2,$$
hence $G'\in \frak
S_{(n-1)/2-r}$ for a positive $r$. When $n=2$, $G^{(\frac12)}\in
\frak S _{2(n-1),0}$ leads to 
$G'\in \frak S_{(n-1)/2,0}$.\qed

\enddemo

We can then conclude:

\proclaim{Theorem 5.15} The eigenvalues of $A_{\chi ,\Sigma _+} ^{-1}-A_\gamma  ^{-1}
$ satisfy:
$$
\mu _j(A_{\chi ,\Sigma _+} ^{-1}-A_\gamma  ^{-1})=\mu _j(\inj_V\gamma _V^{-1}L^{-1}(\gamma _V^{-1})^*\pr_V)=\mu _{j}(P_{2,+}
P_{\nu ,\gamma ,+}P_{2,+} +G),\tag5.35
$$
where $G\in \frak S_{(n-1)/2,0}$.
\endproclaim

\demo{Proof} This follows by inserting the information from Theorem
5.9 in the formula (5.33), using that $P_{2,+}RP_{2,+}\in \frak
S_{(n-1)/2,0}$ by the rules in Section 5.1.\qed
\enddemo

\subhead 5.5 Spectral asymptotics \endsubhead

To find the asymptotic behavior of the $s$-numbers we shall use the following theorem shown in \cite{G11a} (Th.\ 3.3):

\proclaim{Theorem 5.16} {\rm \cite{G11a}} Let $P$ be an operator on $\Sigma $ composed of ${l}$
classical pseudodifferential operators $P_1,\dots, P_{l}$ of negative orders
$-t_1,\dots, -t_{l}$ and  ${l}+1$ functions $b_1,\dots,b_{{l}+1}$ that are
piecewise continuous on $\Sigma $ with possible jumps at $\partial\Sigma  _+$
$$
P=b_1P_1\dots b_{l}P_{l}b_{l+1}.\tag5.36
$$
Let $t=t_1+\dots+t_{l}$. Then $P$ has the spectral behavior:
$$
s_j(P)j^{t/(n-1)} \to c(P)^{t/(n-1)}\text{ for }j\to\infty ,\tag5..37
$$
where 
$$
c(P) 
=\tfrac1{(n-1)(2\pi )^{(n-1)}}\int_{\Sigma }\int_{|\xi '
|=1}
|b_1\dots b_{l+1}p^0_1\dots p^0_l|
^{(n-1)/t}
\,d\omega (\xi ') dx'.
\tag5.38
$$
\endproclaim

Let us also recall  that the principal
symbol of $P_{\nu ,\gamma }$ is $p^0=-|\xi '|^{-1}$. As noted in
Lemma 5.12, the
principal symbol of $P_1=K_\gamma ^*K_\gamma $ is
$\|\tilde k^0\|^2_{L_2}=(2|\xi '|)^{-1}$; that of the squareroot $P_2$ is $\|\tilde
k^0\|_{L_2}=(2|\xi '|)^{-\frac12}$.

Then we can finally show:

\proclaim{Theorem 5.17} Let $\lambda \in \varrho (A_{\chi ,\Sigma _+})\cap
\varrho (A_\gamma )$. The $s$-numbers of $(A_{\chi ,\Sigma _+}-\lambda ) ^{-1}-(A_\gamma -\lambda )^{-1}
$ satisfy the asymptotic formula
$$
s_j((A_{\chi ,\Sigma _+}-\lambda ) ^{-1}-(A_\gamma -\lambda ) ^{-1})j^{2/(n-1)}\to 
C_{0,+}^{2/(n-1)}\text{ for }j\to\infty ,\tag5.39
$$
where 
$$
C_{0,+}= 
\tfrac1{(n-1)(2\pi )^{n-1}}\int_{\Sigma _+}\int_{|\xi
'|=1} (\|\tilde k^0\|_{L_2(\rp)}|p^0|^{1/2})^{n-1}
\,d\omega (\xi ') dx'=c_n\int_{\Sigma _+}1\,dx'.\tag5.40
$$
for a constant $c_n$ depending on $n$ (see {\rm (5.41)} below).
\endproclaim

\demo{Proof} We first treat the case without $\lambda $ (or with
$\lambda =0$), where 
the realizations are positive. Here the $s$-numbers are the positive eigenvalues, and we use (5.35). We can identify $P_{2,+}P_{\nu ,\gamma ,+}P_{2,+} $ with
the operator $1_{\Sigma _+}P_21_{\Sigma _+}P_{\nu ,\gamma }1_{\Sigma _+}P_21_{\Sigma _+}$
in $L_2(\Sigma )$, acting trivially (as 0) on $L_2(\Sigma _-)$. An
application of Theorem 5.16 to this operator gives that
$$
\mu _j(1_{\Sigma _+}P_21_{\Sigma _+}P_{\nu ,\gamma }1_{\Sigma
_+}P_21_{\Sigma _+})j^{2/(n-1)}\to c^{2/(n-1)} \text{ for }j\to\infty ,
$$
where
$$
\aligned
c&=\tfrac1{(n-1)(2\pi )^{(n-1)}}\int_{\Sigma }\int_{|\xi '
|=1}
|1_{\Sigma _+}p^0_2p^0 p^0_2|
^{(n-1)/2}
\,d\omega (\xi ') dx'\\
&=\tfrac1{(n-1)(2\pi )^{(n-1)}}\int_{\Sigma _+}\int_{|\xi '
|=1}
(\|\tilde k^0\|^2|p^0|)
^{(n-1)/2}
\,d\omega (\xi ') dx'=C_{0,+}.
\endaligned
$$
Since $G\in \frak S_{(n-1)/2,0}$, this asymptotic behavior is
preserved under addition of $G$, by Lemma 5.1 $1^\circ$, which implies
the main statement in the theorem for $\lambda =0$. 

Since
$\|\tilde k^0\|^2=(2|\xi '|)^{-1}$, $|p^0|=|\xi '|^{-1}$, the constant
$c_n$ can be calculated as
$$
\aligned
c_n&=\tfrac1{(n-1)(2\pi )^{(n-1)}}\int_{|\xi '
|=1}
2
^{-(n-1)/2}\,d\omega (\xi ') \\
&=\tfrac1{(n-1)(2\pi )^{(n-1)}}
2
^{-(n-1)/2}(n-1)\pi ^{(n-1)/2}\Gamma (1+\tfrac{n-1}2)^{-1}\\
&=(2\pi )^{-(n-1)/2}2^{1-n}\Gamma (1+\tfrac{n-1}2)^{-1}
.
\endaligned\tag5.41
$$

For more general $\lambda \in \varrho (A_{\chi ,\Sigma _+})\cap
\varrho (A_\gamma )$, we use a resolvent identity
 as
in \cite{G11a}:
$$
\multline
(B-\lambda )^{-1}-(B_1-\lambda )^{-1}
=(1+\lambda (B_1-\lambda )^{-1})(B^{-1}-B_1^{-1})(1+\lambda
(B-\lambda )^{-1})\\
=B^{-1}-B_1^{-1} +\lambda (B_1-\lambda
)^{-1}(B^{-1}-B_1^{-1})+(B^{-1}-B_1^{-1}) \lambda  (B-\lambda )^{-1}\\
+ \lambda (B_1-\lambda )^{-1}(B^{-1}-B_1^{-1})\lambda  (B-\lambda )^{-1},
\endmultline \tag5.42$$
valid for $\lambda ,0\in \varrho (B)\cap \varrho
(B_1)$. We apply it to $B=A_{\chi ,\Sigma _+}$ and $B_1=A_\gamma  $ for
$\lambda \in \varrho (\wA)\cap \varrho (A_\gamma  )$. Since $(A_\gamma
-\lambda )^{-1}$ and  $(A_{\chi ,\Sigma _+}-\lambda )^{-1}$ are in $\frak
S_{n/2}$ (cf.\ Corollary 3.2), the
three last terms are in $\frak S_{(n-1)/2 -r}$ with $r>0$. Then we find by
Lemma 5.1 $1^\circ$ that
the main asymptotic estimate of the $s$-numbers is the same as for 
$A_{\chi ,\Sigma _+}^{-1}- A_\gamma ^{-1}$.
\qed
\enddemo

For $n\ge 3$, a generalization of Laptev's result in Theorem 5.2 to
nonstandard $\psi $do's like $\Lambda _+$ and $\Lambda _-$ would allow
an estimate of $s_j(A_{\chi ,\Sigma _+}^{-1}- A_\gamma
^{-1})-C_{0,+}^{2/(n-1)}j^{-2/(n-1)}$ by a lower power of $j$. 

The
methods of \cite{G11} would be useful in an extension of the results
to exterior domains.

\subhead{Acknowledgments}\endsubhead The author is grateful to the referee for encouraging 
inclusion of results on irregular boundaries (with edges). We also
thank Heiko Gimperlein for useful discussions.

\enddocument